\documentclass[12pt,reqno]{amsart}
\usepackage{mathrsfs}
\usepackage{cases}
\usepackage{epic}
\usepackage{amsfonts}
\usepackage{graphicx}
\usepackage{amsmath}
\usepackage{amssymb, upgreek}
\usepackage{bm}
\usepackage{latexsym,todonotes}
\usepackage{pdflscape}
\usepackage[all]{xy}
\usepackage{color}
\usepackage{colordvi}
\usepackage{multicol}
\usepackage{enumitem}
\usepackage[normalem]{ulem}
\newcommand{\arxiv}[1]{\href{http://arxiv.org/abs/#1}{\tt arXiv:\nolinkurl{#1}}}

\usepackage{tikz-cd}
\numberwithin{equation}{section}
\usepackage[linktocpage=true]{hyperref}
\hypersetup{colorlinks,linkcolor=blue,urlcolor=cyan,citecolor=blue}
\usetikzlibrary{arrows.meta}

\newtheorem{theorem}{Theorem}[section]
\newtheorem{proposition}[theorem]{Proposition}
\newtheorem{lemma}[theorem]{Lemma}
\newtheorem{corollary}[theorem]{Corollary}

\theoremstyle{definition}
\newtheorem{definition}[theorem]{Definition}
\newtheorem{example}[theorem]{Example}
\theoremstyle{remark}
\newtheorem{remark}[theorem]{Remark}

\newcommand{\C}{\mathbb C}
\newcommand{\N}{\mathbb N}

\newcommand{\M}{\mathcal M}

\newcommand{\nplus}{\mathfrak n_+}
\newcommand{\g}{\mathfrak g}

\newcommand{\rk}{\rank}
\newcommand{\ad}{\operatorname{ad}}

\newcommand{\Hom}{\operatorname{Hom}}
\newcommand{\Ext}{\operatorname{Ext}}

\newcommand{\diag}{\operatorname{diag}}
\newcommand{\lcm}{\operatorname{lcm}}
\newcommand{\eps}{\varepsilon}

\renewcommand{\Im}{\operatorname{Im}}
\newcommand{\Ker}{\operatorname{Ker}}
\newcommand{\rank}{\operatorname{rank}}
\newcommand{\Rep}{\operatorname{Rep}}
\newcommand{\Gr}{\operatorname{Gr}}
\newcommand{\Fl}{\operatorname{Fl}}
\newcommand{\GL}{\operatorname{GL}}

\theoremstyle{plain}
\newtheorem{introtheorem}{Theorem}

\theoremstyle{plain}

\title[Generalized Preprojective Convolution Algebras]
{Convolution Algebras of Generalized Preprojective Algebras and Kac--Moody Positive Parts}
\author{Fengjun Shi}
\address{Fengjun Shi, Department of Mathematics, Sichuan University, Chengdu 610064, P.~R.~China}
\email{shifengjun@stu.scu.edu.cn}
\thanks{This work was partially supported by the National Natural
Science Foundation of China (Grant No.~12471037).}
\date{}

\begin{document}

\begin{abstract}
We realize the enveloping algebras of positive parts of symmetrizable Kac--Moody Lie algebras as convolution algebras of constructible functions on module varieties of Geiss--Leclerc--Schr\"oer generalized preprojective algebras $\Pi(C,D)$.
We characterize all pairs $(C,D)$ for which such a realization exists.
These realizations generalize Lusztig's classical preprojective construction.
The construction uses a folding embedding and torus actions on flag varieties.
\end{abstract}

    \maketitle
    \setcounter{tocdepth}{1}
    \tableofcontents

\section{Introduction}\label{sec:introduction}
Hall algebras provide representation-theoretic constructions of Lie algebras and their enveloping algebras. For a symmetrizable Kac--Moody Lie algebra $\mathfrak g$, write $\mathbf{U}_v^+(\g)$ for the positive part of the corresponding quantum group. Using representations of hereditary algebras over finite fields, Ringel established a Hall algebra realization of $\mathbf U_v^+(\g)$ in finite type, and Green extended it to the general symmetrizable case \cite{Rin90,Gr95}. The classical specialization at $v=1$ gives $U(\nplus)$, where $\nplus$ is the positive part of $\mathfrak g$.

In the symmetric case, Schofield \cite{Sch91} realized $U(\nplus)$ over $\C$ by replacing finite-field point counts with Euler characteristics of complex varieties. Lusztig \cite{Lus91} reformulated this construction using convolution of constructible functions on quiver representation varieties. He also realized $U(\nplus)$ by convolution on varieties of nilpotent representations of classical preprojective algebras \cite{Lus91,Lus00}.

To extend these geometric constructions to the symmetrizable setting, Geiss, Leclerc, and Schr\"oer \cite{GLS17} introduced the algebras $H(C,D,\Omega)$ and the generalized preprojective algebras $\Pi(C,D)$, where $C$ is a symmetrizable generalized Cartan matrix, $D$ is a positive integral symmetrizer, and $\Omega$ is an acyclic orientation of $C$. Using locally free $H(C,D,\Omega)$-modules, they obtained a convolution realization of $U(\mathfrak n_+(C))$ when $C$ is of Dynkin type \cite{GLS16}. Castillo G\'omez and Geiss \cite{CG25} proved the corresponding result in affine type $\widetilde C_n$ with the minimal symmetrizer. 


For $\Pi(C,D)$, let $\mathcal M(C,D)$ denote the unital convolution subalgebra generated by the characteristic functions $\theta_i=\mathbf{1}_{E_i}$ of the generalized simple $\Pi(C,D)$-modules $E_i$. These generators need not satisfy the Serre relations associated with $C$. In \cite{GLS18}, the authors considered the quotient by the ideal generated by the corresponding Serre elements and proved that it is isomorphic to $U(\nplus(C))$ under their Support Conjecture.

Fei \cite[Theorem~10.1]{Fei26} proved this isomorphism without assuming the Support Conjecture. His proof uses the convolution realization for $H(C,D,\Omega)$, which he established for all symmetrizable $C$ and all positive integral symmetrizers $D$ \cite[Theorem~A.1]{Fei26}. He also gave a counterexample to the conjecture \cite[Theorem~10.2]{Fei26}. In the same work, he realized $U(\nplus(C))$ using monodromy traces on flag varieties for a ramified preprojective species over $\C((t))$.

\subsection{Main results}
We determine when the convolution algebra $\M(C,D)$ is isomorphic to the enveloping algebra of the positive part of a symmetrizable Kac--Moody Lie algebra.

Write $C=(c_{ij})_{i,j\in I}$ and $D=\diag(d_i)_{i\in I}$, and set
\[
    \lambda_i=
    \begin{cases}
        1, & d_i=1,\\
        2, & d_i>1,
    \end{cases}
\]
and define $A=A(C,D)=(a_{ij})_{i,j\in I}$ by
\[
    a_{ii}=2,\qquad
    a_{ij}=\lambda_i c_{ij}\quad(i\ne j).
\]
Then $A$ is a symmetrizable generalized Cartan matrix, with positive integral symmetrizer $\diag(2d_i/\lambda_i)_{i\in I}$. 
Let $e_i$ be the Chevalley generators of $\nplus (A)$, and let $\mathfrak l(C,D)$ be the Lie subalgebra of $\M(C,D)$ generated by $\theta_i$ under the convolution commutator.

For arbitrary $(C,D)$, we prove the $A$-Serre relations in Proposition~\ref{prop:A-Serre}. 
They give a canonical surjective Hopf algebra homomorphism
\[
    \Phi_{C,D}\colon U(\nplus (A)) \twoheadrightarrow \M(C,D), \qquad e_i\mapsto\theta_i.
\]

We consider the following condition, denoted by condition~\eqref{eq:edge-condition}:
\[
    \gcd(d_i,d_j)>1
    \qquad\text{whenever }c_{ij}<0
    \text{ and }d_i,d_j>1.
\]

\begin{introtheorem}[Theorem~\ref{thm:main}]
\label{thm:introduction-main}
If $(C,D)$ satisfies condition~\eqref{eq:edge-condition}, then  
$$\Phi_{C,D}:
    U(\nplus (A))
    \rightarrow
    \M(C,D)$$ 
    is an isomorphism of rank-graded Hopf algebras. Its restriction gives a Lie algebra isomorphism
\[
    \nplus(A)
    \xrightarrow{\;\sim\;}
    \mathfrak l(C,D).
\]
\end{introtheorem}

When $D$ is the identity matrix, we have $A=C$, and Theorem~\ref{thm:introduction-main} recovers Lusztig's classical preprojective realization; see \cite[Section~1.2]{GLS18}. The theorem also applies whenever $C$ is of finite or affine type and $D$ is any positive integral symmetrizer, by Corollary~\ref{cor:finite-affine}. The matrix $A(C,D)$ need not have the same type as $C$.

For pairs $(C,D)$ not satisfying condition~\eqref{eq:edge-condition}, we prove the following result in Section~\ref{app:serre-quotients}.

\begin{introtheorem}[Theorem~\ref{thm:nonrealization}]
\label{thm:introduction-nonrealization}
If $(C,D)$ does not satisfy condition~\eqref{eq:edge-condition}, then
\[
    \mathcal M(C,D)\not\cong U(\mathfrak n_+(B))
\]
as unital $\C$-algebras for any symmetrizable generalized Cartan matrix $B$.
In particular, $\mathfrak l(C,D)$ is not isomorphic to the positive part of any symmetrizable Kac--Moody Lie algebra.
\end{introtheorem}

To obtain realizations for arbitrary $(C,D)$, we construct Serre quotients of $\M(C,D)$ in Section~\ref{app:serre-quotients}. In particular, we recover the full $C$-Serre quotient realization proved by Fei \cite[Theorem~10.1]{Fei26} without assuming the GLS Support Conjecture.

\subsection{Organization}
The paper is organized as follows. 
Section~\ref{sec:preliminaries} recalls generalized preprojective algebras and their convolution algebras, defines the Cartan matrix $A(C,D)$, and states the main theorems.
In Section~\ref{sec:surjectivity}, we prove a splitting lemma for locally free modules and use it to establish the $A$-Serre relations. This gives the canonical homomorphism $U(\mathfrak n_+(A))\to\mathcal M(C,D)$ for arbitrary $(C,D)$.
In Section~\ref{sec:injectivity}, we use a symmetric covering of $A$ and torus actions on flag varieties to prove injectivity under condition ~\eqref{eq:edge-condition}; the flag comparison is related to the acyclic construction in \cite[Appendix~A]{Fei26}.
Section~\ref{app:serre-quotients} proves Theorem~\ref{thm:nonrealization} and constructs Serre quotients for arbitrary $(C,D)$.

\medskip
\noindent\textbf{Acknowledgements.}
The author would like to thank Changjian Fu, Shengfei Geng and Ming Lu for helpful discussions. 
\section{Generalized preprojective algebras and convolution algebras}
\label{sec:preliminaries}
We work over $\C$, and all modules are finite-dimensional left modules.
\subsection{Symmetrizable Cartan matrices and orientations}
    We recall the definitions from \cite{GLS17}. Set $I=\{1,\ldots,n\}$ and $\N=\{0,1,2,\ldots\}$.


    \begin{definition}A matrix $C = (c_{ij}) \in M_n(\mathbb{Z})$ is a symmetrizable (generalized) Cartan matrix if the following conditions are satisfied:
     \begin{itemize}
         \item[($C_1$)] $c_{ii} = 2$ for all $i$;
         \item[($C_2$)] $c_{ij} \leq 0$ for all $i \neq j$;
         \item[($C_3$)] $c_{ij} \neq 0$ if and only if $c_{ji} \neq 0$;
         \item[($C_4$)] There is a diagonal integer matrix $D = \diag(d_1,\cdots,d_n)$ with $d_i \geq 1$ for all $i$ such that $DC$ is symmetric.
     \end{itemize}
     \end{definition}
      The matrix $D$ in $(C_4)$ is called a symmetrizer of $C$. It is minimal if $d_1 +\cdots +d_n$ is minimal. For each pair with $c_{ij}<0$, set
     \begin{align*}
         g_{ij}:=|\gcd(c_{ij},c_{ji})|, \quad f_{ij}:=|c_{ij}|/g_{ij}, \quad k_{ij}:= \gcd(d_{i},d_{j}).
     \end{align*}
     Since $d_if_{ij}=d_jf_{ji}$ and $\gcd(f_{ij},f_{ji})=1$, we have
     \[
     d_i=k_{ij}f_{ji}, \qquad d_j=k_{ij}f_{ij}.
     \]

     An orientation of $C$ is a subset $\Omega \subseteq I\times I$ satisfying the following conditions:
     \begin{enumerate}[label=(\roman*)]
         \item $\{(i, j), (j, i)\} \cap \Omega \neq \emptyset $ if and only if $c_{ij} < 0$;
         \item For each sequence $((i_1, i_2), (i_2, i_3), \cdots , (i_t, i_{t+1}))$ with $t \geq 1$ and $(i_s, i_{s+1}) \in \Omega$ for all $1 \leq s \leq t$, we have $i_1 \neq i_{t+1}$.
     \end{enumerate}

     For an orientation $\Omega$ of $C$, let $Q=Q(C,\Omega)$ be the quiver with vertex set $Q_0=I$ and arrow set
     \[
         Q_1= 
    \left\{
    \alpha_{ij}^{(g)}:j\to i \;\middle|\; (i,j)\in\Omega,\ 1\leq g\leq g_{ij} \right\}
    \cup 
    \left\{
    \varepsilon_i:i\to i \;\middle|\; i\in I \right\}.
     \]
     Thus, $Q$ has $g_{ij}$ arrows from $j$ to $i$ for each $(i,j)\in\Omega$, and one loop $\varepsilon_i$ at each vertex. The quiver obtained by deleting all loops $\varepsilon_i$ is acyclic.

     Set
     \[
         \Omega^{\mathrm{op}} =\{(j,i)\mid(i,j)\in\Omega\}, \qquad \overline{\Omega}=\Omega\cup\Omega^{\mathrm{op}}.
     \]

    Let $\overline{Q}$ be the quiver obtained from $Q$ by adding a reverse arrow $\alpha_{ji}^{(g)}:i\to j$ for each non-loop arrow $\alpha_{ij}^{(g)}:j\to i$. The loops $\varepsilon_i$ are retained without adding new copies. For $(i,j)\in\overline{\Omega}$, define
    \[
    s_{ij}= \begin{cases}
    1, & (i,j)\in\Omega,\\
    -1, & (i,j)\in\Omega^{\mathrm{op}}.
    \end{cases}
    \]

\subsection{Generalized preprojective algebras}

    \begin{definition}\label{def:H}
    The algebra $H=H(C,D,\Omega)$ is the quotient
    \[
    H(C,D,\Omega)=\mathbb{C}Q/\mathcal{I}_H,
    \]
    where $\mathcal{I}_H$ is the two-sided ideal generated by the following relations:
    \begin{enumerate}
    \item For each $i\in I$, the nilpotency relation
    \[
    \varepsilon_i^{d_i}=0.
    \]
    \item For each $(i,j)\in\Omega$ and
    $1\leq g\leq g_{ij}$, the commutativity relation
    \[
    \varepsilon_i^{f_{ji}}\alpha_{ij}^{(g)}
    =
    \alpha_{ij}^{(g)}\varepsilon_j^{f_{ij}}.
    \]
    \end{enumerate}
    \end{definition}


    \begin{definition}\label{def:Pi}
    The generalized preprojective algebra $\Pi=\Pi(C,D,\Omega)$ is the quotient
    \[
    \Pi(C,D,\Omega)
    =
    \mathbb{C}\overline{Q}/\mathcal{I}_{\Pi},
   \]
   where $\mathcal{I}_{\Pi}$ is the two-sided ideal generated by the following relations:
    \begin{enumerate}
    \item For each $i\in I$, the nilpotency relation
    \[
    \varepsilon_i^{d_i}=0.
    \]
    \item For each $(i,j)\in\overline{\Omega}$ and
    $1\leq g\leq g_{ij}$, the commutativity relation
    \[
    \varepsilon_i^{f_{ji}}\alpha_{ij}^{(g)}
    =
    \alpha_{ij}^{(g)}\varepsilon_j^{f_{ij}}.
    \]
    \item For each $i\in I$, the mesh relation
    \[
    \sum_{\substack{j\in I\\c_{ij}<0}}
    s_{ij}
    \sum_{g=1}^{g_{ij}}
    \sum_{p=0}^{f_{ji}-1}
    \varepsilon_i^p
    \alpha_{ij}^{(g)}\alpha_{ji}^{(g)}
    \varepsilon_i^{f_{ji}-1-p}
    =0.
    \]
    \end{enumerate}
    \end{definition}

    The isomorphism class of $\Pi(C,D,\Omega)$ is independent of $\Omega$, so we write $\Pi(C,D)$. If $C$ is symmetric and $D=I$, the loops vanish and $\Pi(C,D)$ is the ordinary preprojective algebra.

\subsection{The quadratic form}
    We associate with $(C,D)$ the quadratic form $q_{C,D}:\mathbb{Z}^n\to\mathbb{Z}$ defined by
    \[
    q_{C,D}(\mathbf{v})
    =
    \frac{1}{2}\mathbf{v}^{T}DC\mathbf{v}
    =
    \sum_{i=1}^{n}d_i v_i^2 - \sum_{i<j}d_i|c_{ij}|v_i v_j,
    \qquad
    \mathbf{v}=(v_1,\ldots,v_n).
    \]
    The matrix $C$ is of finite (Dynkin) type if and only if the extension of $q_{C,D}$ to $\mathbb{R}^n$ is positive definite. The connected finite types are
    \[
    A_n,\ B_n,\ C_n,\ D_n,\ E_6,\ E_7,\ E_8,\ F_4,\ G_2.
    \]
    If $C$ is connected and of affine type, this form is positive semidefinite but not positive definite.
    
\subsection{Locally free and generalized simple modules}
   Put $H_i=\C[\eps_i]/(\eps_i^{d_i})$. For a $\Pi(C,D)$-module $M$, write $M_i$ for its vector space at vertex $i$, regarded as an $H_i$-module via the action of $\varepsilon_i$. A module $M$ is \emph{locally free} if $M_i$ is free over $H_i$ for every $i$, and its rank vector is
   \[
   \rk M=(\rk_{H_i}M_i)_{i\in I}.
   \]
   The generalized simple module $E_i$ is supported at $i$ with component $H_i$. It is indecomposable; it is an ordinary simple module only when $d_i=1$.

   A module $M$ is called \emph{$E$-filtered} if it admits
   a finite filtration
   \[
   0=M^0\subset M^1\subset\cdots\subset M^m=M
   \]
   such that
   \[
   M^a/M^{a-1}\cong E_{i_a}, \qquad 1\leq a\leq m,
   \]
   for some $i_1,\ldots,i_m\in I$.

   Every $E$-filtered module is locally free. Indeed, at each vertex $j\in I$, the filtration induces short exact sequences
   \[
   0\longrightarrow (M^{a-1})_j \longrightarrow (M^a)_j \longrightarrow (E_{i_a})_j \longrightarrow 0.
   \]
   The quotient $(E_{i_a})_j$ is either $H_j$ or zero, so each sequence splits as a sequence of $H_j$-modules. It follows inductively that $M_j$ is free over $H_j$.

\subsection{The convolution algebra}
   Following  \cite[Section~4.1] {GLS16} and \cite[Section~2.5] {GLS18}, we recall the construction of the convolution algebra of constructible functions associated with $\Pi:=\Pi(C,D)$. 

   For an arrow $a$ of $\overline Q$, let $s(a)$ and $t(a)$ denote its source and target. For a dimension vector $\mathbf v=(v_i)_{i\in I}\in\mathbb N^I$, let $\Rep_{\mathbf v}(\Pi)$ be the affine variety of $\Pi$-module structures on the vertex spaces $(\mathbb C^{v_i})_{i\in I}$. Concretely, its points are tuples
   \[
   M=(M(a))_{a\in\overline Q_1} \in \prod_{a\in\overline Q_1} \Hom_{\mathbb C} \bigl(\mathbb C^{v_{s(a)}}, \mathbb C^{v_{t(a)}}\bigr)
   \]
   satisfying the defining relations of $\Pi$.

   The group
   \[
   G_{\mathbf v}
   =
   \prod_{i\in I}\GL_{v_i}(\mathbb C)
   \]
   acts on $\Rep_{\mathbf v}(\Pi)$ by changes of bases:
   \[
   (g\cdot M)(a)
   =
   g_{t(a)}M(a)g_{s(a)}^{-1}, \qquad g=(g_i)_{i\in I}\in G_{\mathbf v}.
   \]
   Its orbits correspond bijectively to isomorphism classes of $\Pi$-modules with dimension vector $\mathbf v$.

   A constructible function on a complex algebraic variety is a complex-valued function taking finitely many values and having constructible fibers. Let $\mathcal F_{\mathbf v}(\Pi)$ be the $\mathbb C$-vector space of constructible functions
   \[
   f:\Rep_{\mathbf v}(\Pi)\longrightarrow \mathbb C
   \]
   that are constant on $G_{\mathbf v}$-orbits; equivalently,
   \[
   f(g\cdot M)=f(M) \qquad \text{for all }g\in G_{\mathbf v} \text{ and }M\in\Rep_{\mathbf v}(\Pi).
   \]
   Set
   \[ \mathcal F(\Pi)
   =
   \bigoplus_{\mathbf v\in\mathbb N^I} \mathcal F_{\mathbf v}(\Pi).
   \]
   We regard these functions as functions on isomorphism classes of modules, extending each $f\in\mathcal F_{\mathbf v}(\Pi)$ by zero to modules whose dimension vectors differ from $\mathbf v$.

   Let $\chi_c$ denote the compactly supported Euler characteristic. For a constructible function $\varphi:X\to\mathbb C$, define its Euler integral by
   \[ \int_X \varphi\,d\chi_c
   =
   \sum_{c\in\Im(\varphi)} c\,\chi_c\bigl(\varphi^{-1}(c)\bigr).
   \]
   This sum is finite since $\varphi$ takes finitely many values.

   For a finite-dimensional $\Pi$-module $M$, let
   \[
   \Gr_{\mathbf v}^{\Pi}(M)
   =
   \{U\subseteq M: U \text{ is a }\Pi\text{-submodule and }
  \underline{\dim}U=\mathbf v\}.
   \]
   This is a closed subvariety of $\prod_{i\in I}\Gr_{v_i}(M_i)$, and hence is projective.

   For $f\in\mathcal F_{\mathbf v}(\Pi)$ and $g\in\mathcal F_{\mathbf w}(\Pi)$, define
   \begin{equation}\label{eq:Pi-convolution}
   (f*g)(M)
   =
   \int_{\Gr_{\mathbf v}^{\Pi}(M)} f(U)\,g(M/U)\,d\chi_c.
   \end{equation}
   $(\mathcal F(\Pi),*)$ is an $\mathbb N^I$-graded associative algebra. Its identity is the characteristic function $\mathbf 1_0$ of the zero module.
   
   For each $i\in I$, let $\theta_i$ be the characteristic function of the isomorphism class of the generalized simple module $E_i$:
   \[
   \theta_i(M)
   =
   \begin{cases}
   1,& M\cong E_i,\\
   0,& \text{otherwise}.
   \end{cases}
   \]
   Since $\underline{\dim}E_i=d_i\mathbf e_i$, we have $\theta_i\in\mathcal F_{d_i\mathbf e_i}(\Pi)$, where $(\mathbf e_i)_{i\in I}$ is the standard basis of $\mathbb Z^I$.

   \begin{definition}
   Let $\mathcal M(C,D)$ be the unital subalgebra of $(\mathcal F(\Pi),*)$ generated by the functions $\theta_i$:
   \[
   \mathcal M(C,D)
   =
   \mathbb C\langle\theta_i:i\in I\rangle_*.
   \]
   \end{definition}

   This is the algebra denoted by $\widetilde{\mathcal M}(\Pi)$ in \cite[Section~2.5]{GLS18}.

   We equip $\mathcal M(C,D)$ with the rank grading $\deg\theta_i=\mathbf e_i$. Equivalently, set
   \[
   \mathcal M_{\mathbf r}
   =
   \mathcal M(C,D)\cap\mathcal F_{D\mathbf r}(\Pi),
   \qquad
   \mathcal M(C,D)
   =
   \bigoplus_{\mathbf r\in\mathbb N^I}
   \mathcal M_{\mathbf r}.
   \]
   Each $\mathcal M_{\mathbf r}$ is finite-dimensional: it is spanned by the finitely many convolution words in which $\theta_i$ occurs $r_i$ times. Also, $\mathcal M_{\mathbf 0}=\mathbb C\mathbf 1_0$.

   For a sequence $\mathbf i=(i_1,\ldots,i_m)$, let $\Fl_{\mathbf i}^{E}(M)$ be the constructible set of flags of $\Pi$-submodules
   \[
   0=M^0\subset M^1\subset\cdots\subset M^m=M
   \]
   such that
   \[
   M^s/M^{s-1}\cong E_{i_s} \qquad (1\le s\le m).
   \]
   Iterating the convolution formula gives
   \begin{equation}\label{eq:Pi-flag-evaluation}
   (\theta_{i_1}*\cdots*\theta_{i_m})(M)
   =
   \chi_c\bigl(\Fl_{\mathbf i}^{E}(M)\bigr).
   \end{equation}

   Consequently, every element of $\mathcal M(C,D)$ vanishes on modules that are not $E$-filtered. The following is the analogue of \cite[Lemma~4.2]{GLS16} for $\Pi$. It is also recorded in \cite[Section~2.5]{GLS18}.
   \begin{lemma}\label{lem:locally-free-support}
   Let $X$ be a finite-dimensional $\Pi$-module. If $X$ is not locally free, then $f(X)=0$ for every $f\in\mathcal M(C,D)$.
   \end{lemma}
   
   We recall the following Hopf algebra structure from \cite[Section~2.5]{GLS18}; see also \cite[Propositions~4.5 and~4.7]{GLS16} for the corresponding construction for $H$.
   
   \begin{proposition}\label{prop:Hopf-structure}
   The convolution algebra $\mathcal M(C,D)$ is a connected cocommutative Hopf algebra. Its coproduct
   \[
   \Delta:\mathcal M(C,D)\longrightarrow \mathcal M(C,D)\otimes\mathcal M(C,D)
   \]
   is given by
   \[
   (\Delta f)(X,Y)=f(X\oplus Y).
   \]
   In particular,
   \[ \Delta(\theta_i) = \theta_i\otimes\mathbf1_0 + \mathbf1_0\otimes\theta_i.
   \]
   The counit is $\epsilon(f)=f(0)$, and the antipode satisfies $S(\theta_i)=-\theta_i$.
\end{proposition}
   Since $\Delta$ is an algebra homomorphism, for any sequence $\mathbf i=(i_1,\ldots,i_m)$ we have
   \[
   \Delta(\theta_{i_1}*\cdots*\theta_{i_m})
   =
   \prod_{a=1}^{m} \bigl(
   \theta_{i_a}\otimes\mathbf1_0 + \mathbf1_0\otimes\theta_{i_a} \bigr).
   \]

   \begin{definition}\label{def:primitive-element}
   Let $B$ be a Hopf algebra with coproduct $\Delta$ and unit $1_B$. An element $x\in B$ is called \emph{primitive} if
   \[
   \Delta(x)=x\otimes 1_B+1_B\otimes x.
   \]
   We denote the vector space of primitive elements by
   \[
   P(B)= \{x\in B:\Delta(x)=x\otimes 1_B+1_B\otimes x\}.
   \]
   \end{definition}
   
   The space $P(B)$ is a Lie subalgebra of $B$ under the commutator $[x,y]=xy-yx$. 
   Since each $\theta_i$ is primitive by Proposition~\ref{prop:Hopf-structure}, every iterated convolution commutator of the $\theta_i$ is also primitive.

   For the convolution algebra, primitive elements have the following interpretation; see also \cite[Section~2.5]{GLS18}.

   \begin{lemma}\label{lem:primitive-support}
   A function $f\in\mathcal M(C,D)$ is primitive if and only if
   \[
   f(0)=0 \qquad\text{and}\qquad f(X\oplus Y)=0 \quad\text{for all nonzero $\Pi$-modules $X,Y$}.
   \]
   Equivalently, a primitive function can take a nonzero value only on a nonzero indecomposable module.
\end{lemma}

\begin{proof}
    By Proposition~\ref{prop:Hopf-structure}, a function $f\in\mathcal M(C,D)$ is primitive if and only if
\[
   f(X\oplus Y)
   =
   f(X)\mathbf1_0(Y)+\mathbf1_0(X)f(Y)
\]
   for all finite-dimensional $\Pi$-modules $X,Y$. If $f$ is primitive, taking $X=Y=0$ gives $f(0)=0$, and taking $X,Y\ne0$ gives $f(X\oplus Y)=0$. Conversely, these two vanishing conditions imply the displayed identity: both sides vanish when $X,Y\ne0$, and the identity is immediate when one of the modules is zero.
\end{proof}
   As recalled in \cite[Section~2.5]{GLS18}, the inclusion
\[
  P(\mathcal M(C,D))\hookrightarrow\mathcal M(C,D)
\]
  induces a canonical Hopf algebra isomorphism
\[
  U\bigl(P(\mathcal M(C,D))\bigr) \xrightarrow{\sim} \mathcal M(C,D).
\]
Let $\mathfrak l(C,D)$ denote the Lie subalgebra of $\mathcal M(C,D)$ generated by the $\theta_i$ under the convolution commutator.
\begin{lemma}\label{lem:enveloping-realization}
    The inclusion $\mathfrak l(C,D)\hookrightarrow\mathcal M(C,D)$ induces a canonical Hopf algebra isomorphism
\[
    U(\mathfrak l(C,D))\xrightarrow{\sim}\mathcal M(C,D).
\]
\end{lemma}
\begin{proof}
    Since each $\theta_i$ is primitive, we have $\mathfrak l(C,D)\subseteq P(\mathcal M(C,D))$. By the PBW theorem, the induced map
    \[
    U(\mathfrak l(C,D))\longrightarrow\mathcal M(C,D)
\]
is injective. Its image contains all the algebra generators $\theta_i$, so it is surjective. 
\end{proof}

\subsection{The modified Cartan matrix and the main theorems}
\label{subsec:main-theorem}
Let $C=(c_{ij})_{i,j\in I}$ be a symmetrizable generalized Cartan matrix, and let $D=\operatorname{diag}(d_i)$ be a positive integral symmetrizer. Set
\[
    \lambda_i=
    \begin{cases}
        1, & d_i=1,\\
        2, & d_i>1,
    \end{cases}
\]
and define $A=A(C,D)=(a_{ij})_{i,j\in I}$ by
\begin{equation}\label{eq:modified-cartan}
    a_{ii}=2,\qquad
    a_{ij}=\lambda_i c_{ij}\quad(i\ne j).
\end{equation}
Then $A$ is a symmetrizable generalized Cartan matrix
with positive integral symmetrizer
\[
    D_A=\diag(2d_i/\lambda_i)_{i\in I}.
\]
We consider the condition
\begin{equation}\tag{E}\label{eq:edge-condition}
    \gcd(d_i,d_j)>1
    \qquad\text{whenever }c_{ij}<0
    \text{ and }d_i,d_j>1.
\end{equation}

Let $e_i$, $i\in I$, be the Chevalley generators of $\mathfrak n_+(A)$. We can now state the realization theorem.

\begin{theorem}\label{thm:main}
Assume that $(C,D)$ satisfies
condition~\eqref{eq:edge-condition}.
Then the assignment $e_i\mapsto\theta_i$ extends uniquely
to an isomorphism of rank-graded Hopf algebras
\[
    \Phi_{C,D}\colon
    U(\mathfrak n_+(A))
    \xrightarrow{\;\sim\;}
    \mathcal M(C,D).
\]
Its restriction gives a Lie algebra isomorphism
\[
    \mathfrak n_+(A)\xrightarrow{\;\sim\;} \mathfrak l(C,D).
\]
\end{theorem}

In Section~\ref{sec:surjectivity}, we prove the Serre relations associated with $A$ and establish surjectivity for arbitrary $(C,D)$, without assuming condition~\eqref{eq:edge-condition}. In Section~\ref{sec:injectivity}, we prove injectivity under this condition, completing the proof of Theorem~\ref{thm:main}.

\begin{corollary}[Finite and affine types] \label{cor:finite-affine}
Let $C$ be of finite or affine type, and let $D$ be any positive integral symmetrizer.  We have a canonical rank-graded Hopf algebra isomorphism
\[
\begin{aligned}
    U(\mathfrak n_+(A)) &\xrightarrow{\sim} \mathcal M(C,D),\\
    e_i &\longmapsto \theta_i,
\end{aligned}
\]
where $A=A(C,D)$. In particular, the Lie subalgebra generated by the $\theta_i$ under the convolution commutator is canonically isomorphic to $\mathfrak n_+(A)$.
\end{corollary}

\begin{proof}
Since $C$ is of finite or affine type, the symmetric matrix $DC$ is positive semidefinite. Its two-by-two principal minors therefore
give
\[
    c_{ij}c_{ji}\leq 4
    \qquad (i\neq j).
\]
It is then easy to check that condition~\eqref{eq:edge-condition} holds. The assertions now follow from Theorem~\ref{thm:main}.
\end{proof}

For any symmetrizable generalized Cartan matrix $C$ and any positive integral symmetrizer $D'$, the pair $(C,mD')$ satisfies condition~\eqref{eq:edge-condition} for every integer $m\ge 2$, since all diagonal entries of $mD'$ are divisible by $m$. Thus Theorem~\ref{thm:main} applies to $(C,mD')$.

For pairs $(C,D)$ not satisfying condition~\eqref{eq:edge-condition}, we prove the following theorem in Section~\ref{app:serre-quotients}.

\begin{theorem}\label{thm:nonrealization}
   If $(C,D)$ does not satisfy condition~\eqref{eq:edge-condition}, then
\[
    \mathcal M(C,D)\not\cong U(\mathfrak n_+(B))
\]
as unital $\mathbb C$-algebras for any symmetrizable generalized Cartan matrix $B$. In particular, $\mathfrak l(C,D)$ is not isomorphic to the positive part of any symmetrizable Kac--Moody Lie algebra. 
\end{theorem}

\section{Serre relations and the canonical homomorphism} \label{sec:surjectivity}
In this section, we prove the Serre relations associated with $A=A(C,D)$ and obtain the canonical homomorphism to $\M(C,D)$. We do not assume condition~\eqref{eq:edge-condition}.
We use the incoming and outgoing maps from \cite[Sections~5.1--5.3]{GLS17}. 

For $c_{ij}<0$, let ${}_iH_j$ denote the $(H_i,H_j)$-bimodule of arrows from $j$ to $i$. It is free of rank $|c_{ij}|$ as a left $H_i$-module. For a $\Pi$-module $M$, set
\[
    \widetilde M_i
    =
    \bigoplus_{j:\,c_{ij}<0}
    {}_iH_j\otimes_{H_j}M_j.
\]
Following \cite[Section~5.3]{GLS17}, we consider the $H_i$-module homomorphisms
\[
    M_{i,\mathrm{in}}\colon\widetilde M_i\longrightarrow M_i,
    \qquad
    M_{i,\mathrm{out}}\colon M_i\longrightarrow\widetilde M_i.
\]
By \cite[Proposition~5.2]{GLS17}, the preprojective relations are equivalent to
\[
    M_{i,\mathrm{in}}\cdot M_{i,\mathrm{out}}=0
    \qquad(i\in I).
\]

\subsection{Splitting lemma}\label{subsec:splitting}
We adapt the splitting argument in the proof of \cite[Lemma~6.3]{GLS18} to locally free $\Pi$-modules.
\begin{lemma}\label{lem:general splitting}
   Let $X$ be a locally free $\Pi$-module and let $i\in I$. Let $m$ and $n$ be the minimal numbers of $H_i$-generators of $\Im X_{i,\mathrm{in}}$ and $\Im X_{i,\mathrm{out}}$. If $\rank_{H_i}X_i>m+n$, then $E_i$ is a direct summand of $X$.
\end{lemma}

\begin{proof}
    Write $\alpha=X_{i,\mathrm{in}}$ and $\beta=X_{i,\mathrm{out}}$. Thus
\[
    \widetilde X_i \xrightarrow{\alpha} X_i
    \xrightarrow{\beta} \widetilde X_i,
    \qquad \alpha\beta=0.
\]
We put the matrix of $\alpha$ in Smith normal form. Since $\Im \alpha$ is minimally generated by $m$ elements, this form has exactly $m$ nonzero diagonal entries. Thus we obtain a decomposition into free $H_i$-modules
\[
    X_i=P\oplus Z,\qquad
    P\cong H_i^m,\qquad
    \Im \alpha\subseteq P.
\]
Similarly, applying Smith normal form to $\beta$, we obtain a decomposition into free $H_i$-modules
\[
    \widetilde X_i=T\oplus T',\qquad
    T\cong H_i^n,\qquad
    \Im \beta\subseteq T.
\]
Now consider the restriction $\beta|_Z:Z\to T$. Since
\[
    \rank_{H_i}Z
    =\rank_{H_i}X_i-m>n=\rank_{H_i}T,
\]
we may choose a direct summand $K\cong H_i$ of $Z$ such that $\beta(K)=0$. Since $\Im \alpha\subseteq P$, $E_i$ is a direct summand of $X$.
\end{proof}

\begin{corollary}\label{lem:splitting}
Let $i\ne j$ with $c_{ij}<0$, and let $X$ be a locally free $\Pi$-module of rank $r\mathbf e_i+\mathbf e_j$. If $r>\lambda_i|c_{ij}|$, then $E_i$ is a direct summand of $X$.
\end{corollary}

\begin{proof}
    Put $N=|c_{ij}|$. Since $X$ is locally free, we have
\begin{align*}
    X_i\cong H_{i}^r,\qquad X_j\cong H_j,\qquad
    \widetilde X_i
    ={}_iH_j\otimes_{H_j}X_j\cong H_i^N.
\end{align*}
Write $\alpha=X_{i,\mathrm{in}}$ and $\beta=X_{i,\mathrm{out}}$. Let $m,n$ be as in Lemma~\ref{lem:general splitting}.

If $d_i=1$, then $H_i=\mathbb C$. Since
$\operatorname{Im}\beta\subseteq\ker\alpha$, we have
\[
    m+n
    =\operatorname{rank}\alpha+\operatorname{rank}\beta
    \leq N.
\]
If $d_i>1$, we have
\[
    H_i^N \xrightarrow{\alpha} H_i^r \xrightarrow{\beta} H_i^N.
\]
Thus $m \leq N$ and $n \leq N$.

In both cases, $m+n\leq\lambda_i N<r$. Thus $E_i$ is a direct summand of $X$ by Lemma~\ref{lem:general splitting}.
\end{proof}

\subsection{The Serre relations}\label{subsec:Serre}
We now apply Corollary~\ref{lem:splitting} to prove the Serre relations associated with $A(C,D)$. For $f,g\in\mathcal M(C,D)$, write
\[
   \ad(f)(g)=[f,g]=f*g-g*f.
\]

\begin{proposition}\label{prop:A-Serre}
The generators $\theta_i$, $i\in I$, satisfy
\begin{equation}\label{eq:A-Serre}
(\ad\theta_i)^{1-a_{ij}}(\theta_j)=0 \qquad (i\ne j).
\end{equation}
\end{proposition}

\begin{proof}
    Fix $i \neq j$ and put
\[
   r=1-a_{ij}=1+\lambda_i|c_{ij}|, \qquad \Theta_{ij}=(\ad\theta_i)^r(\theta_j).
\]
$\Theta_{ij}$ is primitive and homogeneous of rank degree $r\mathbf e_i+\mathbf e_j$. By Lemma~\ref{lem:locally-free-support} and Lemma~\ref{lem:primitive-support}, it can be nonzero only on locally free indecomposable modules of this rank.

Let $X$ be a locally free module with rank $r\mathbf e_i+\mathbf e_j$. If $c_{ij}=0$, then $r=1$ and $X \cong E_i \oplus E_j$. If $c_{ij}<0$, Corollary ~\ref{lem:splitting} gives $X\cong E_i\oplus Y$. Thus each such $X$ is decomposable and $\Theta_{ij}=0$.
\end{proof}

\begin{corollary}\label{cor:canonical-surjection}
For every symmetrizable generalized Cartan matrix $C$ and every positive integral symmetrizer $D$, the assignment $e_i\mapsto\theta_i$ extends uniquely to a surjective rank-graded Hopf algebra homomorphism
\begin{equation}\label{eq:canonical-surjection}
\Phi_{C,D}\colon U\bigl(\nplus(A(C,D))\bigr) \twoheadrightarrow \mathcal M(C,D).
\end{equation}
\end{corollary}

\begin{proof}
 The Serre presentation and Proposition~\ref{prop:A-Serre} give the unique algebra homomorphism with $e_i\mapsto\theta_i$. It is surjective because $\theta_i$ generate $\M(C,D)$ and it preserves the grading and the Hopf structure because the generators on both sides are primitive of degree $\mathbf e_i$.
\end{proof}

\section{Injectivity of the canonical homomorphism} \label{sec:injectivity}

Put $A=A(C,D)$. We first construct an auxiliary quiver $Q$ for an arbitrary pair $(C,D)$. Under condition~\eqref{eq:edge-condition}, 
each representation $V$ of $Q$ gives a locally free $\Pi(C,D)$-module $F(V)$. To prove that $\Phi_{C,D}$ is injective, we use a torus action to compare the $E$-flags of $F(V)$ with flags of $Q$-subrepresentations of $V$. This comparison is closely related to the argument in \cite[Appendix~A.2]{Fei26}. In the construction of $F(V)$, we choose weight shifts for the arrow maps and use condition~\eqref{eq:edge-condition} to verify the mesh relations.

\subsection{Quiver convolution algebras and symmetric coverings}
Let $Q$ be a finite quiver without loops, with vertex set $Q_0$; oriented cycles and multiple arrows are allowed. Let $B=(b_{ij})_{i,j\in Q_0}$ be its associated symmetric generalized Cartan matrix:
\[
  b_{ii}=2,\qquad b_{ij} = -\#\{\text{arrows joining $i$ and $j$ in either direction}\} \quad(i\ne j).
\]

For $\mathbf d\in\mathbb N^{Q_0}$, let $\mathcal F_{\mathbf d}(Q)$ be the space of
isomorphism-invariant $\mathbb C$-valued constructible functions on $\Rep_{\mathbf d}(Q)$. The space
\[
  \mathcal F(Q) = \bigoplus_{\mathbf d\in\mathbb N^{Q_0}} \mathcal F_{\mathbf d}(Q)
\]
  is equipped with the convolution product
\[
  (f*g)(V) = \int_{W\subseteq V} f(W)g(V/W)\,d\chi_c,
\]
where $W$ runs through the subrepresentations of $V$.

For each $i\in Q_0$, let $S_i$ be the one-dimensional representation supported at $i$, and put $s_i=\mathbf1_{S_i}$. We denote by
\[
  \mathcal M(Q) = \mathbb C\langle s_i:i\in Q_0\rangle_* \subseteq\mathcal F(Q)
\]
the unital convolution subalgebra generated by the $s_i$.

We recall the classical constructible-function realization of the positive enveloping algebra, originating in the work of Schofield and Lusztig \cite{Sch91,Lus91}. For an explicit statement for quivers without loops, see \cite[Theorem~3.1]{JS12}.

\begin{theorem}\label{thm:inj-hall}
The assignment $e_i\mapsto s_i$ extends uniquely
to a graded algebra isomorphism
\[
\eta_Q\colon U(\nplus(B))
\xrightarrow{\sim}\mathcal M(Q).
\]
\end{theorem}


   For a sequence $\mathbf j=(j_1,\ldots,j_r)$ of vertices of $Q$ and a finite-dimensional representation $V$, let $\Fl^{Q}_{\mathbf j}(V)$ denote the variety of flags of subrepresentations
\[
   0=V^0\subset V^1\subset\cdots\subset V^r=V
\]
  such that
\[
  V^k/V^{k-1}\simeq S_{j_k} \qquad (1\leq k\leq r).
\]
  With our convolution convention, we have
\begin{equation}\label{eq:inj-Q-flags}
  (s_{j_1}*\cdots*s_{j_r})(V) = \chi_c\bigl(\Fl^{Q}_{\mathbf j}(V)\bigr).
\end{equation}
For integers $m>0$ and $w$, write $[w]_{m}$ for the least nonnegative residue of $w$ modulo $m$. The following folding lemma is due to Hubery \cite[Theorem~7.1.5]{Hub02}; see also \cite[Proposition~4.4]{Dem11}.

\begin{lemma}\label{lem:inj-cover}
Let $A=(a_{ij})_{i,j\in I}$ be a symmetrizable generalized Cartan matrix, and let $B$ be a symmetric generalized Cartan matrix indexed by
\[
\{(i,a):i\in I,\ 0\leq a<m_i\}.
\]
Here each $m_i$ is a positive integer. Assume that the permutation $\sigma(i,a)=(i,[a+1]_{m_i})$ preserves $B$, that no two distinct vertices over the same $i$ are adjacent, and that
\begin{equation}\label{eq:inj-columnsum}
\sum_{a=0}^{m_i-1}b_{(i,a),(j,c)}
=
a_{ij}
\qquad
(i,j\in I,\ 0\leq c<m_j).
\end{equation}
Then $e_i\longmapsto \sum_{a=0}^{m_i-1}e_{(i,a)}$ defines an injective Lie algebra homomorphism
\[
\iota\colon
\nplus(A)\hookrightarrow\nplus(B).
\]
The induced homomorphism $U(\iota)$ is also injective.
\end{lemma}

\begin{proof}
The cyclic group $\langle\sigma\rangle$ acts admissibly on the diagram of $B$, and \eqref{eq:inj-columnsum} identifies its folding Cartan matrix with $A$. The Lie algebra embedding is therefore given by \cite[Proposition~4.4]{Dem11}, and PBW gives the assertion about enveloping algebras.
\end{proof}

\begin{remark}
When $B$ comes from a quiver, $\sigma$ need only preserve the underlying graph, not its arrow directions.
\end{remark}

\subsection{Construction of the auxiliary quiver and the folding embedding}\label{subsec:auxiliary-quiver}
  We construct a finite quiver $Q$ without loops whose associated symmetric generalized Cartan matrix $B$ folds to $A=A(C,D)$. The construction and the embedding below do not require condition~\eqref{eq:edge-condition}.

\subsubsection{Vertices}
Set
\begin{equation}\label{eq:inj-weightsizes}
  L=\lcm\bigl(2,d_i:i\in I\bigr), \qquad
  m_i=
\begin{cases}
  L/2,& d_i=1,\\
  L/d_i,& d_i>1.
\end{cases}
\end{equation}
  Define
\[
    Q_0=\{(i,a):i\in I,\ 0\le a<m_i\},
    \qquad
    I_i=\{(i,a):0\le a<m_i\}.
\]

\subsubsection{Arrows}
For each unordered adjacent pair $\{i,j\}$ of $C$, apply exactly one of the following rules. In each case, $1\le\ell\le g_{ij}$.

\par\smallskip
\noindent\textbf{Case 1: $d_i=d_j=1$.}
Here $m_i=m_j=L/2$. Choose and fix a direction, say $j\to i$. For each $0\leq w<L/2$
introduce an arrow
\[
\alpha_{ij,w}^{(\ell)}:
(j,w)\longrightarrow(i,w).
\]

\par\smallskip
\noindent\textbf{Case 2: exactly one of $d_i,d_j$ equals $1$.}
After interchanging $i$ and $j$ if necessary, assume that $d_i>1$ and $d_j=1$. Here $m_i=L/d_i$ and $m_j=L/2$. For each $0\le w<L/2$, introduce the arrows
\begin{equation}\label{eq:inj-mixed-arrows}
\begin{aligned}
\beta_{ij,w}^{(\ell)}:
(i,[w]_{m_i})&\longrightarrow(j,w),\\
\alpha_{ij,w}^{(\ell)}:
(j,w)&\longrightarrow(i,[w+L/2]_{m_i}).
\end{aligned}
\end{equation}

\par\smallskip
\noindent\textbf{Case 3: $d_i,d_j>1$.}
Here $m_i=L/d_i$ and $m_j=L/d_j$. Recall that $k_{ij}=\gcd(d_i,d_j)$, and put
\[
t=L/k_{ij}=\lcm(m_i,m_j).
\]
For each $0\le w<t$, introduce the arrows
\begin{equation}\label{eq:auxiliary-long-arrows}
\begin{aligned}
    \alpha^{(\ell)}_{ij,w}&\colon
    (j,[w]_{m_j})\longrightarrow(i,[w]_{m_i}),\\
    \beta^{(\ell)}_{ij,w}&\colon
    (i,[w]_{m_i})\longrightarrow(j,[w]_{m_j}).
\end{aligned}
\end{equation}
There are no other arrows.

\subsubsection{The folding embedding}

Let
\[
    B=\bigl(b_{(i,a),(j,b)}\bigr)
\]
be the symmetric generalized Cartan matrix associated with $Q$.
For each $i\in I$, set
\[
    \sigma_i=\sum_{a=0}^{m_i-1}s_{(i,a)}
    \in\mathcal M(Q).
\]

\begin{theorem}\label{thm:auxiliary-embedding}
The assignment
\[
    e_i\longmapsto\sum_{a=0}^{m_i-1}e_{(i,a)}
\]
extends to an injective Lie algebra homomorphism
\[
    \iota\colon
    \mathfrak n_+(A)\hookrightarrow\mathfrak n_+(B).
\]
Consequently, the assignment $e_i\mapsto\sigma_i$
extends uniquely to an injective algebra homomorphism
\[
    J=\eta_Q\circ U(\iota)\colon
    U(\mathfrak n_+(A))\hookrightarrow\mathcal M(Q).
\]
\end{theorem}

\begin{proof}
    No arrow joins two vertices in the same fiber $I_i$. We check that $\sigma(i,a)=(i,[a+1]_{m_i})$ preserves $B$ and that
    \begin{align}
        \sum_{a=0}^{m_i-1}b_{(i,a),(j,b)}=a_{ij}.
    \end{align}
The identity is immediate when $i=j$ or $c_{ij}=0$. We therefore assume that $c_{ij}<0$ and consider the following three cases.
\par\smallskip
\noindent\textbf{Case 1: $d_i=d_j=1$.}
Here $m_i=m_j=L/2$. 
For each $0\le w<L/2$, there are $g_{ij}$ edges between $(i,w)$ and $(j,w)$. $\sigma$ preserves these edges and
\[
    \sum_{a=0}^{m_i-1}b_{(i,a),(j,b)}
    =-g_{ij}=c_{ij}=a_{ij}.
\]
The same argument applies with $i$ and $j$ interchanged.

\par\smallskip
\noindent\textbf{Case 2: exactly one of $d_i,d_j$ equals $1$.} After interchanging $i$ and $j$ if necessary, assume that $d_i>1$ and $d_j=1$. Then $m_i=L/d_i$ and $m_j=L/2$, and
\[
    c_{ij}=-g_{ij},
    \qquad
    c_{ji}=-d_i g_{ij}.
\]

For each $0\le b<m_j$, the vertex $(j,b)$ is joined to $(i,[b]_{m_i})$ and $(i,[b+m_j]_{m_i})$ by $g_{ij}$ edges each, counted with multiplicity. After forgetting the arrow directions, $\sigma$ sends each edge in either family at $b$ to the corresponding edge in the same family at $b+1$, provided $b<m_j-1$. At $b=m_j-1$, it sends the edge underlying $\beta^{(\ell)}_{ij,m_j-1}$ to that underlying $\alpha^{(\ell)}_{ij,0}$, and the edge underlying $\alpha^{(\ell)}_{ij,m_j-1}$ to that underlying $\beta^{(\ell)}_{ij,0}$, since $2m_j=L$ is divisible by $m_i$. Thus $\sigma$ preserves the edge multiplicities between the two fibers.

Each vertex $(j,b)$ is incident with $2g_{ij}$ edges joining it to $I_i$, so
\[
    \sum_{a=0}^{m_i-1}b_{(i,a),(j,b)}
    =-2g_{ij}
    =2c_{ij}
    =a_{ij}.
\]
The indices $b$ and $b+m_j$, as $b$ ranges from $0$ to $m_j-1$, together run through $0,\ldots,L-1$. Each residue modulo $m_i$ therefore occurs $L/m_i=d_i$ times. Hence
\[
    \sum_{b=0}^{m_j-1}b_{(j,b),(i,a)}
    =-d_i g_{ij}
    =c_{ji}
    =a_{ji}.
\]

\par\smallskip
\noindent\textbf{Case 3: $d_i,d_j>1$.}
Put
\[
    t=\frac{L}{k_{ij}}
      =\lcm(m_i,m_j).
\]
Each index $0\le w<t$ contributes $2g_{ij}$ edges between
$(i,[w]_{m_i})$ and $(j,[w]_{m_j})$. Since $m_i$ and $m_j$ divide $t$, $\sigma$ corresponds to replacing $w$ by $[w+1]_t$ and preserves the edge multiplicities.

For each fixed $0\le b<m_j$, exactly $t/m_j$ values of $w$ satisfy $[w]_{m_j}=b$. Since
\[
    \frac{t}{m_j}
    =\frac{d_j}{k_{ij}}
    =f_{ij},
\]
we obtain
\[
    \sum_{a=0}^{m_i-1}b_{(i,a),(j,b)}
    =-2g_{ij}\frac{t}{m_j}
    =-2g_{ij}f_{ij}
    =2c_{ij}
    =a_{ij}.
\]
Similarly, $t/m_i=d_i/k_{ij}=f_{ji}$, and hence
\[
    \sum_{b=0}^{m_j-1}b_{(j,b),(i,a)}
    =-2g_{ij}\frac{t}{m_i}
    =-2g_{ij}f_{ji}
    =2c_{ji}
    =a_{ji}.
\]
The result now follows from Lemma~\ref{lem:inj-cover} and Theorem~\ref{thm:inj-hall}.
\end{proof}

\begin{remark}
After forgetting the arrow directions, the construction of $Q$ admits a uniform description. For each unordered adjacent pair $\{i,j\}$, set
\[
    t_{ij}=\lcm(m_i,m_j),
    \qquad
    \rho_{ij}=\gcd(|a_{ij}|,|a_{ji}|).
\]
For each $0\le w<t_{ij}$, join $(i,[w]_{m_i})$ and $(j,[w]_{m_j})$ by $\rho_{ij}$ edges. Then $\sigma$ preserves the edge multiplicities, and
\[
    \rho_{ij}\frac{t_{ij}}{m_j}=|a_{ij}|,
    \qquad
    \rho_{ij}\frac{t_{ij}}{m_i}=|a_{ji}|.
\]
This gives the folding relations. We use the chosen arrow directions to construct locally free $\Pi(C,D)$-modules from $Q$-representations in the next subsection.
\end{remark}

\begin{example}\label{ex:auxiliary-quiver}
Let
\[
    D=\diag(1,1,2,4),
    \qquad
    C=
    \begin{pmatrix}
        2&-1&0&0\\
        -1&2&-2&0\\
        0&-1&2&-2\\
        0&0&-1&2
    \end{pmatrix}.
\]
Then $DC$ is symmetric, and
\[
    A(C,D)=
    \begin{pmatrix}
        2&-1&0&0\\
        -1&2&-2&0\\
        0&-2&2&-4\\
        0&0&-2&2
    \end{pmatrix}.
\]
We have $g_{12}=g_{23}=g_{34}=1$ and
\[
    L=4,
    \qquad
    (m_1,m_2,m_3,m_4)=(2,2,2,1).
\]
Thus
\[
    Q_0=
    \{(i,a):1\le i\le3,\ a=0,1\}
    \sqcup\{(4,0)\}.
\]
Choosing the direction $2\to1$ in Case~1, the resulting quiver is
\[
    (1,0)
    \xleftarrow{\alpha_{12,0}}
    (2,0)
    \overset{\alpha_{32,0}}
      {\underset{\beta_{32,0}}{\rightleftarrows}}
    (3,0)
    \overset{\beta_{34,0}}
      {\underset{\alpha_{34,0}}{\rightleftarrows}}
    (4,0)
    \overset{\alpha_{34,1}}
      {\underset{\beta_{34,1}}{\rightleftarrows}}
    (3,1)
    \overset{\beta_{32,1}}
      {\underset{\alpha_{32,1}}{\rightleftarrows}}
    (2,1)
    \xrightarrow{\alpha_{12,1}}
    (1,1).
\]
\end{example}

\subsection{From quiver representations to
generalized preprojective modules}
Assume condition~\eqref{eq:edge-condition}. Let $V$ be a finite-dimensional representation of the auxiliary quiver $Q$. Write $V_{i,a}$ for its space at $(i,a)$ and $V(\rho)$ for its linear map along an arrow $\rho$. We construct a locally free $\Pi(C,D)$-module $F(V)$
by specifying its vertex spaces and arrow actions.

\subsubsection{Vertex spaces and weights}
For each $i\in I$, define
\begin{equation}\label{eq:inj-testvertices}
F(V)_i
=
\bigoplus_{a=0}^{m_i-1}
H_i\otimes_{\mathbb C}V_{i,a}.
\end{equation}
The loop $\varepsilon_i$ acts by multiplication
on the $H_i$-factor.

For $0\le a<m_i$ and $0\le p<d_i$, assign weight $a+m_ip$ to the subspace
\[
    F_i[a+m_ip]
    :=
    \mathbb C\varepsilon_i^p
    \otimes_{\mathbb C}V_{i,a}.
\]
The admissible weights are
\[
    \begin{cases}
        0\le w<L/2, & d_i=1,\\
        0\le w<L,   & d_i>1.
    \end{cases}
\]
For all other integers $w$, set $F_i[w]=0$. For every admissible weight $w=a+m_ip$, we use
the $\mathbb C$-linear isomorphism
\[
    F_i[w]\xrightarrow{\;\sim\;}V_{i,[w]_{m_i}},
    \qquad
    \varepsilon_i^p\otimes v\longmapsto v.
\]
The loop action satisfies
\[
\varepsilon_i F_i[w]
\subseteq F_i[w+m_i].
\]

\subsubsection{Arrow actions}

For each unordered adjacent pair $\{i,j\}$,
use the same labeling as in the construction
of $Q$.
For $1\leq\ell\leq g_{ij}$, we define maps
\[
a_\ell:F(V)_j\longrightarrow F(V)_i,
\qquad
b_\ell:F(V)_i\longrightarrow F(V)_j
\]
as the actions of the arrows
$\alpha_{ij}^{(\ell)}$ and
$\alpha_{ji}^{(\ell)}$ of $\Pi(C,D)$,
respectively.

The maps on weight spaces below are defined via the canonical $\mathbb C$-linear identifications
\[
    F_i[w]\cong V_{i,[w]_{m_i}}.
\]
All unlisted weight blocks are zero.

\par\smallskip
\noindent\textbf{Case 1: $d_i=d_j=1$.}
Retain the direction $j\to i$ chosen in the construction of $Q$.
For $0\leq w<L/2$, define
\[
a_\ell\big|_{F_j[w]}:
F_j[w]
\xrightarrow{\;V(\alpha_{ij,w}^{(\ell)})\;}
F_i[w],
\]
and set $b_\ell=0$.

\par\smallskip
\noindent\textbf{Case 2: $d_i>1$ and $d_j=1$.}
For $0\leq w<L/2$, define
\begin{equation}\label{eq:inj-mixedplacement}
\begin{aligned}
b_\ell\big|_{F_i[w]}:\quad
F_i[w]
&\xrightarrow{\;V(\beta_{ij,w}^{(\ell)})\;}
F_j[w],\\
a_\ell\big|_{F_j[w]}:\quad
F_j[w]
&\xrightarrow{\;V(\alpha_{ij,w}^{(\ell)})\;}
F_i[w+L/2].
\end{aligned}
\end{equation}
In particular, $b_\ell$ is zero on all weight spaces of weight at least $L/2$.

\par\smallskip
\noindent\textbf{Case 3: $d_i,d_j>1$.}
Put
\[
t=\frac{L}{k_{ij}},
\qquad
s=L-t,
\]
where $k_{ij}=\gcd(d_i,d_j)$. For $0\leq w<t$, define
\begin{equation}\label{eq:inj-thickplacement}
\begin{aligned}
a_\ell\big|_{F_j[w]}:\quad
F_j[w]
&\xrightarrow{\;V(\alpha_{ij,w}^{(\ell)})\;}
F_i[w+s],\\
b_\ell\big|_{F_i[w]}:\quad
F_i[w]
&\xrightarrow{\;V(\beta_{ij,w}^{(\ell)})\;}
F_j[w+s].
\end{aligned}
\end{equation}
Both maps are zero on source weight spaces of weight at least $t$.

Since $t=\lcm(m_i,m_j)$ and both
$m_i$ and $m_j$ divide $L$, the shift $s=L-t$
is divisible by both $m_i$ and $m_j$.
Consequently,
\[
    [w+s]_{m_i}=[w]_{m_i},
    \qquad
    [w+s]_{m_j}=[w]_{m_j}.
\]
Thus the source and target weight spaces identify with the source and target spaces of the corresponding arrow maps of $V$.

\begin{lemma}\label{lem:relations}
Assume condition~\eqref{eq:edge-condition}. For every finite-dimensional representation $V$ of $Q$, the operators defined above endow
\[
    F(V)=\bigoplus_{i\in I}F(V)_i
\]
with the structure of a locally free $\Pi(C,D)$-module.
Moreover,
\[
    \rank_{H_i}F(V)_i
    =
    \sum_{a=0}^{m_i-1}\dim_{\mathbb C}V_{i,a}
    \qquad(i\in I).
\]
\end{lemma}

\begin{proof}
The relations $\varepsilon_i^{d_i}=0$ hold by the definition of $H_i$.
We prove the stronger assertion that both sides of every commutativity relation and each individual summand in the mesh relations act as zero on $F(V)$. We check this for each adjacent pair $\{i,j\}$ and each $1\le\ell\le g_{ij}$.
\par\smallskip
\noindent\textbf{Case 1: $d_i=d_j=1$.}
The arrows of $Q$ between $I_i$ and $I_j$ all have the chosen direction $j\to i$. For the reverse arrow $\alpha_{ji}^{(\ell)}$ of $\Pi(C,D)$, we set $b_\ell=0$.
Since $d_i=d_j=1$, both loops act as zero, so both sides of the commutativity relations vanish. Moreover, $b_\ell=0$ implies
\[
    a_\ell b_\ell=b_\ell a_\ell=0.
\]

\par\smallskip
\noindent\textbf{Case 2: $d_i>1$ and $d_j=1$.}
Here $f_{ij}=1$ and $f_{ji}=d_i$.
Since $\varepsilon_i^{d_i}=0$ and $\varepsilon_j=0$, both sides of the commutativity relations vanish:
\[
    \varepsilon_i^{d_i}a_\ell
    =0=a_\ell\varepsilon_j,
    \qquad
    \varepsilon_jb_\ell
    =0=b_\ell\varepsilon_i^{d_i}.
\]

The image of $a_\ell$ lies in weights
at least $L/2$, on which $b_\ell$ is zero.
Thus
\[
b_\ell a_\ell=0,
\]
which gives the contribution to the mesh relation at $j$.

At $i$, the corresponding mesh terms are
\[
\varepsilon_i^p a_\ell b_\ell
\varepsilon_i^{d_i-1-p},
\qquad 0\leq p<d_i.
\]
The map $a_\ell b_\ell$ raises weight by
$L/2$, while each loop raises weight by
$m_i$. Hence each displayed term raises
weight by
\[
\frac L2+(d_i-1)m_i
=
\frac{3L}{2}-m_i
\geq L,
\]
where we used $d_i\geq2$ and $m_i=L/d_i$. Since $F_i[w]=0$ outside $0\le w<L$, every such mesh term is zero.

\par\smallskip
\noindent\textbf{Case 3: $d_i,d_j>1$.}
Recall that
\[
    t=\frac{L}{k_{ij}},
    \qquad
    s=L-t,
\]
and
\[
    d_i=k_{ij}f_{ji},
    \qquad
    d_j=k_{ij}f_{ij}.
\]
Hence
\[
    f_{ji}m_i=f_{ij}m_j=t,
\]
so the loop powers $\varepsilon_i^{f_{ji}}$ and $\varepsilon_j^{f_{ij}}$ both raise weight by $t$.

Both $a_\ell$ and $b_\ell$ raise weight by $s$, while $\varepsilon_i^{f_{ji}}$ and $\varepsilon_j^{f_{ij}}$ raise weight by $t$. Hence each side of the commutativity relations raises weight by $s+t=L$. Since all admissible weights at $i$ and $j$ lie between $0$ and $L-1$, we obtain
\[
    \varepsilon_i^{f_{ji}}a_\ell
    =a_\ell\varepsilon_j^{f_{ij}}=0,
    \qquad
    \varepsilon_j^{f_{ij}}b_\ell
    =b_\ell\varepsilon_i^{f_{ji}}=0.
\]

At vertex $i$, consider a mesh term
\[
    \varepsilon_i^p a_\ell b_\ell
    \varepsilon_i^{f_{ji}-1-p},
    \qquad 0\le p<f_{ji}.
\]
Each arrow map raises weight by $s$, while each application of $\varepsilon_i$ raises
weight by $m_i$. Hence this term raises weight by
\[
    2s+(f_{ji}-1)m_i
    =
    2(L-t)+t-m_i
    =
    2L-t-m_i.
\]
Under condition~\eqref{eq:edge-condition}, we have $k_{ij}\ge2$, so $t=L/k_{ij}\le L/2$. Moreover, $m_i\le t$ since $t=f_{ji}m_i$. Consequently,
\[
    2L-t-m_i\ge 2L-2t\ge L.
\]
Since $F_i[w]=0$ outside $0\le w<L$, the displayed mesh term is zero. Interchanging $i$ and $j$ gives the same conclusion for every mesh term at $j$.

Summing these zero terms with the prescribed signs proves the mesh relations at every vertex. Together with the loop and commutativity relations, this shows that the operators define a $\Pi(C,D)$-module structure on $F(V)$.

Finally, each summand $H_i\otimes_{\mathbb C}V_{i,a}$ is free over $H_i$ of rank $\dim_{\mathbb C}V_{i,a}$. Thus $F(V)$ is locally free, with
\[
    \rank_{H_i}F(V)_i
    =
    \sum_{a=0}^{m_i-1}\dim_{\mathbb C}V_{i,a}
    \qquad(i\in I).
\]
\end{proof}

\begin{example}\label{ex:auxiliary-module}
Continue with Example~\ref{ex:auxiliary-quiver}.
Take $V_{i,a}=\mathbb C$ at every vertex of $Q$ and assign
arbitrary scalars $x_a,y_a,z_a,u_a,v_a\in\mathbb C$,
for $a=0,1$, to its arrows as follows:
\[
\begin{tikzcd}[column sep=1.5em]
    \underset{(1,0)}{\mathbb C}
    & \underset{(2,0)}{\mathbb C}
        \arrow[l,"x_0"']
        \arrow[r,shift left=.6ex,"y_0"]
    & \underset{(3,0)}{\mathbb C}
        \arrow[l,shift left=.6ex,"z_0"]
        \arrow[r,shift left=.6ex,"v_0"]
    & \underset{(4,0)}{\mathbb C}
        \arrow[l,shift left=.6ex,"u_0"]
        \arrow[r,shift left=.6ex,"u_1"]
    & \underset{(3,1)}{\mathbb C}
        \arrow[l,shift left=.6ex,"v_1"]
        \arrow[r,shift left=.6ex,"z_1"]
    & \underset{(2,1)}{\mathbb C}
        \arrow[l,shift left=.6ex,"y_1"]
        \arrow[r,"x_1"]
    & \underset{(1,1)}{\mathbb C}
\end{tikzcd}
\]
Here each arrow label denotes multiplication by the
indicated scalar, and the label beneath each space
specifies its vertex.

Using the weight-ordered bases
\[
\begin{aligned}
    \mathcal B_1&=(e_{1,0},e_{1,1}),\\
    \mathcal B_2&=(e_{2,0},e_{2,1}),\\
    \mathcal B_3&=(e_{3,0},e_{3,1},
        \varepsilon_3e_{3,0},\varepsilon_3e_{3,1}),\\
    \mathcal B_4&=(e_4,\varepsilon_4e_4,
        \varepsilon_4^2e_4,\varepsilon_4^3e_4),
\end{aligned}
\]
where $e_{i,a}=1\otimes1$ and $e_4=e_{4,0}$,
we obtain the following representation of $\Pi(C,D)$:
\[
\begin{tikzcd}[column sep=3em]
    \underset{1}{\mathbb C^2}
        \arrow[r,shift right=.6ex,"0"']
    & \underset{2}{\mathbb C^2}
        \arrow[l,shift right=.6ex,"X_{12}"']
        \arrow[r,shift left=.6ex,"X_{32}"]
    & \underset{3}{\mathbb C^4}
        \arrow[l,shift left=.6ex,"X_{23}"]
        \arrow[r,shift left=.6ex,"X_{43}"]
        \arrow[loop above,"N_3"]
    & \underset{4}{\mathbb C^4}
        \arrow[l,shift left=.6ex,"X_{34}"]
        \arrow[loop above,"N_4"]
\end{tikzcd}
\]
Here the loop operators are
\[
    N_3=
    \begin{pmatrix}
        0&0&0&0\\
        0&0&0&0\\
        1&0&0&0\\
        0&1&0&0
    \end{pmatrix},
    \qquad
    N_4=
    \begin{pmatrix}
        0&0&0&0\\
        1&0&0&0\\
        0&1&0&0\\
        0&0&1&0
    \end{pmatrix}.
\]
The loops at vertices $1$ and $2$ act as zero and are omitted from the diagram. With respect to the chosen bases and using column vectors, the arrow matrices are
\[
\begin{aligned}
    X_{12}
    &=
    \begin{pmatrix}
        x_0&0\\
        0&x_1
    \end{pmatrix},
    &
    X_{21}
    &=
    \begin{pmatrix}
        0&0\\
        0&0
    \end{pmatrix},
    \\[3mm]
    X_{32}
    &=
    \begin{pmatrix}
        0&0\\
        0&0\\
        y_0&0\\
        0&y_1
    \end{pmatrix},
    &
    X_{23}
    &=
    \begin{pmatrix}
        z_0&0&0&0\\
        0&z_1&0&0
    \end{pmatrix},
    \\[3mm]
    X_{34}
    &=
    \begin{pmatrix}
        0&0&0&0\\
        0&0&0&0\\
        u_0&0&0&0\\
        0&u_1&0&0
    \end{pmatrix},
    &
    X_{43}
    &=
    \begin{pmatrix}
        0&0&0&0\\
        0&0&0&0\\
        v_0&0&0&0\\
        0&v_1&0&0
    \end{pmatrix}.
\end{aligned}
\]
Here $X_{ij}$ represents the action of the arrow $\alpha_{ij}\colon j\to i$. The module $F(V)$ has rank vector $(2,2,2,1)$ and complex dimension vector $(2,2,4,4)$.
\end{example}

\subsection{Torus action and flag comparison}
\label{subsec:flag-comparison}
Let $T=\mathbb C^\times$ act on $F(V)$ through the weight decomposition by
\[
    g_i(t)v=t^w v
    \qquad
    (v\in F_i[w],\ t\in T).
\]
Write $g(t)=\bigoplus_{i\in I}g_i(t)$. 

A $\mathbb C$-linear map $\rho\colon F(V)_j\to F(V)_i$ is called homogeneous of degree $\delta\in\mathbb Z$ if
\[
    \rho(F_j[w])\subseteq F_i[w+\delta]
    \qquad\text{for every }w\in\mathbb Z.
\]

\begin{lemma}\label{lem:torus-action}
For every sequence $\mathbf i$, the weight action
on $F(V)$ induces an algebraic action of
$T=\mathbb C^\times$ on
$\mathrm{Fl}_{\mathbf i}^{E}(F(V))$, given by
\[
    t\cdot(U^\bullet)=g(t)U^\bullet.
\]
\end{lemma}

\begin{proof}
    The loop $\varepsilon_i$ has degree $m_i$, and the arrow degrees are
\[
    (\deg a_\ell,\deg b_\ell)=
    \begin{cases}
        (0,0),
            & d_i=d_j=1,\\
        (L/2,0),
            & d_i>1,\ d_j=1,\\
        (L-L/k_{ij},L-L/k_{ij}),
            & d_i,d_j>1.
    \end{cases}
\]
For every loop or arrow operator $\rho\colon F(V)_j\to F(V)_i$ of degree $\delta_\rho$, we have
\[
    g_i(t)\rho g_j(t)^{-1}
    =t^{\delta_\rho}\rho.
\]
Let $U=(U_i)_{i\in I}$ be a
$\Pi(C,D)$-submodule of $F(V)$. Then
\[
    \rho\bigl(g_j(t)U_j\bigr)
    =
    t^{-\delta_\rho}g_i(t)\rho(U_j)
    \subseteq g_i(t)U_i.
\]
Thus $g(t)U$ is again a submodule.

Let $U'\subseteq U$ be submodules with $U/U'\cong E_i$. The quotient $g(t)U/g(t)U'$ is supported at vertex $i$, so all arrows between distinct vertices act as zero. Under the linear identification induced by $g(t)$, its loop action is $t^{-m_i}$ times the original one, since
\[
    \varepsilon_i g(t)
    =t^{-m_i}g(t) \varepsilon_i.
\]
Rescaling the loop action by a nonzero scalar preserves the isomorphism class of $E_i$, so
\[
    g(t)U/g(t)U'\cong E_i.
\]
Applying this to each successive quotient shows that $g(t)$ preserves $\Fl_{\mathbf i}^{E}(F(V))$.
\end{proof}

Applying the Euler-characteristic fixed-point formula \cite[Corollary~2]{BB73} to the action in Lemma~\ref{lem:torus-action}, we obtain
\begin{equation}\label{eq:euler-fixed}
    \chi_c\bigl(\mathrm{Fl}_{\mathbf i}^{E}(F(V))\bigr)
    =
    \chi_c\bigl(\mathrm{Fl}_{\mathbf i}^{E}(F(V))^T\bigr).
\end{equation}
The fixed locus consists precisely of the flags whose terms are $T$-stable submodules. To describe these flags, we first identify the $T$-stable locally free submodules of $F(V)$.

\begin{lemma}\label{lem:free}
Let $U_i\subseteq F(V)_i$ be a $T$-stable $H_i$-submodule which is free over $H_i$. Then there are unique subspaces $W_{i,a}\subseteq V_{i,a}$ such that
\[
    U_i=
    \bigoplus_{a=0}^{m_i-1}
    H_i\otimes_{\mathbb C}W_{i,a}.
\]
Conversely, every submodule of this form is $T$-stable and free, and its quotient in $F(V)_i$ is free.
\end{lemma}

\begin{proof}
Since $U_i$ is $T$-stable, it decomposes into weight spaces:
\[
    U_i=\bigoplus_w\bigl(U_i\cap F_i[w]\bigr).
\]
The weights $a+pm_i$ are pairwise distinct, so there are unique subspaces
$W_{i,a,p}\subseteq V_{i,a}$ such that
\[
    U_i=
    \bigoplus_{a=0}^{m_i-1}
    \bigoplus_{p=0}^{d_i-1}
    \mathbb C\varepsilon_i^p
    \otimes_{\mathbb C}W_{i,a,p}.
\]
Stability under $\varepsilon_i$ gives
\[
    W_{i,a,0}\subseteq W_{i,a,1}
    \subseteq\cdots\subseteq W_{i,a,d_i-1}
    \qquad(0\le a<m_i).
\]

Let $r=\rank_{H_i}U_i$. Since $U_i$ is free over $H_i$, we have
\[
    r=\dim_{\mathbb C}
    \bigl(\varepsilon_i^{d_i-1}U_i\bigr)=\dim_{\mathbb C}\bigl( \bigoplus_{a=0}^{m_i-1}
    \mathbb C\varepsilon_i^{d_i-1}
    \otimes_{\mathbb C}W_{i,a,0} \bigr)=\sum_{a=0}^{m_i-1}
    \dim_{\mathbb C}W_{i,a,0}.
\]
Thus
\[
\begin{aligned}
    d_i r
    &=\dim_{\mathbb C}U_i\\
    &=\sum_{a=0}^{m_i-1}
      \sum_{p=0}^{d_i-1}
      \dim_{\mathbb C}W_{i,a,p}\\
    &\ge
      d_i\sum_{a=0}^{m_i-1}
      \dim_{\mathbb C}W_{i,a,0}
      =d_i r.
\end{aligned}
\]
Hence we have $W_{i,a,p}=W_{i,a,0}$ for all $a,p$. Taking $W_{i,a}=W_{i,a,0}$ proves the decomposition and uniqueness.

Conversely, every submodule of the stated form is free over $H_i$ and $T$-stable. Moreover,
\[
    F(V)_i/U_i
    \cong
    \bigoplus_{a=0}^{m_i-1}
    H_i\otimes_{\mathbb C}(V_{i,a}/W_{i,a}),
\]
so the quotient is also free over $H_i$.
\end{proof}

We now describe the $T$-fixed flags in terms of flags of $Q$-subrepresentations.

\begin{lemma}\label{lem:flag-comparison}
Let $V$ be a finite-dimensional representation of $Q$, and let $\mathbf i=(i_1,\ldots,i_r)$. There is a natural isomorphism of varieties
\[
    \mathrm{Fl}_{\mathbf i}^{E}(F(V))^T
    \cong
    \bigsqcup_{\substack{0\le a_k<m_{i_k}\\1\le k\le r}}
    \mathrm{Fl}^{Q}_{((i_1,a_1),\ldots,(i_r,a_r))}(V).
\]
Consequently,
\[
    \chi_c\bigl(\mathrm{Fl}_{\mathbf i}^{E}(F(V))\bigr)
    =
    \sum_{\substack{0\le a_k<m_{i_k}\\1\le k\le r}}
    \chi_c\bigl(
        \mathrm{Fl}^{Q}_{((i_1,a_1),\ldots,(i_r,a_r))}(V)
    \bigr).
\]
\end{lemma}

\begin{proof}
Let $U^\bullet=\{0=U^0\subset U^1\subset\cdots\subset U^r=F(V)\} \in \Fl_{\mathbf i}^{E}(F(V))^T$. Each $U_i^k$ is $T$-stable and free over $H_i$. Lemma~\ref{lem:free} therefore gives unique subspaces $W_{i,a}^k\subseteq V_{i,a}$ such that
\[
    U_i^k=
    \bigoplus_{a=0}^{m_i-1}
    H_i\otimes_{\mathbb C}W_{i,a}^k.
\]
The inclusions $U^{k-1}\subseteq U^k$ imply $W_{i,a}^{k-1}\subseteq W_{i,a}^k$. By the construction of $F(V)$, the submodule condition on $U^k$ implies that $W^k=(W_{i,a}^k)$ is a $Q$-subrepresentation of $V$. Thus $W^\bullet=\{ 0=W^0\subset W^1\subset\cdots\subset W^r=V \}$ is a flag of $Q$-subrepresentations. Moreover,
\begin{align*}
    U_i^k/U_i^{k-1}
    \cong
    \bigoplus_{a=0}^{m_i-1}
    H_i\otimes_{\mathbb C}
    (W_{i,a}^k/W_{i,a}^{k-1}).
\end{align*}
Since $U^k/U^{k-1}\cong E_{i_k}$, there is a unique $0\le a_k<m_{i_k}$ for which $W_{i_k,a_k}^k/W_{i_k,a_k}^{k-1}$ is one-dimensional, and all other quotient spaces are zero. As $Q$ has no loops, we obtain
\begin{align*}
    W^k/W^{k-1}\cong S_{(i_k,a_k)}.
\end{align*}
Conversely, let $W^\bullet \in \Fl^{Q}_{((i_1,a_1),\ldots,(i_r,a_r))}(V)$. Let $U_i^k= \bigoplus_{a=0}^{m_i-1}H_i\otimes_{\mathbb C}W_{i,a}^k$.
By Lemma~\ref{lem:free}, each $U_i^k$ is a $T$-stable free $H_i$-submodule of $F(V)_i$. Since $W^k$ is a $Q$-subrepresentation, the construction of $F(V)$ shows that $U^k=(U_i^k)_{i\in I}$ is a $\Pi(C,D)$-submodule. The inclusions in $W^\bullet$ give
\begin{align*}
    0=U^0\subset U^1\subset\cdots\subset U^r=F(V).
\end{align*}
Since $W^k/W^{k-1}\cong S_{(i_k,a_k)}$, as $H_i$-modules we have
\[
\begin{aligned}
    U_i^k/U_i^{k-1}
    &\cong
    \bigoplus_{a=0}^{m_i-1}
    H_i\otimes_{\mathbb C}
    (W_{i,a}^k/W_{i,a}^{k-1})\\
    &\cong
    \begin{cases}
        H_i,&i=i_k,\\
        0,&i\ne i_k.
    \end{cases}
\end{aligned}
\]
Thus $U^k/U^{k-1}\cong E_{i_k}$, so $U^\bullet$ is a $T$-fixed $E$-flag. The two constructions are inverse by the uniqueness in Lemma~\ref{lem:free}.

It remains to show that both maps are regular.  For a fixed sequence $(a_1,\ldots,a_r)$, the dimensions of all spaces $W_{i,a}^k$ are fixed. The map $W^\bullet\mapsto U^\bullet$ is regular, since it is given by tensor products with the fixed spaces $H_i$ and direct sums.
Its inverse recovers $W_{i,a}^k$ from
\[
    U_i^k\cap F_i[a]
    =
    \mathbb C1\otimes_{\mathbb C}W_{i,a}^k.
\]
Since these intersections have fixed dimensions, the inverse map is also regular.

The dimensions $\dim_{\mathbb C}(U_i^k\cap F_i[w])$
are locally constant on the fixed locus of the
ambient flag variety.
Hence the loci indexed by $(a_1,\ldots,a_r)$
are open and closed in
$\mathrm{Fl}_{\mathbf i}^{E}(F(V))^T$. This proves the isomorphism of varieties.

Finally, the Euler-characteristic identity follows from \eqref{eq:euler-fixed} and additivity of $\chi_c$.
\end{proof}

We can now complete the proof of the main theorem.

\begin{proof}[Proof of Theorem~\ref{thm:main}]
By Corollary~\ref{cor:canonical-surjection}, $\Phi_{C,D}$ is a surjective rank-graded Hopf algebra homomorphism. Let $J:U(\mathfrak n_+(A))\hookrightarrow \mathcal M(Q)$ be the injection in Theorem~\ref{thm:auxiliary-embedding}. For $\mathbf i=(i_1,\ldots,i_r)$ and every $Q$-representation $V$,  Lemma~\ref{lem:flag-comparison} gives
\[
\begin{aligned}
    &(\theta_{i_1}*\cdots*\theta_{i_r})(F(V))\\
    &\qquad=
    \chi_c\bigl(\mathrm{Fl}_{\mathbf i}^{E}(F(V))\bigr)\\
    &\qquad=
    \sum_{\substack{0\le a_k<m_{i_k}\\1\le k\le r}}
    \chi_c\bigl(
        \mathrm{Fl}^{Q}_{((i_1,a_1),\ldots,(i_r,a_r))}(V)
    \bigr)\\
    &\qquad=
    \sum_{\substack{0\le a_k<m_{i_k}\\1\le k\le r}}
    (s_{(i_1,a_1)}*\cdots*s_{(i_r,a_r)})(V)\\
    &\qquad=
    (\sigma_{i_1}*\cdots*\sigma_{i_r})(V).
\end{aligned}
\]
By linearity, we obtain
\begin{equation}\label{eq:evaluation-comparison}
    \Phi_{C,D}(u)(F(V))=J(u)(V)
\end{equation}
for every $u\in U(\mathfrak n_+(A))$ and every finite-dimensional representation $V$ of $Q$.

Now let $u\in\ker\Phi_{C,D}$. Equation~\eqref{eq:evaluation-comparison} gives
\[
    J(u)(V)=\Phi_{C,D}(u)(F(V))=0
\]
for every finite-dimensional representation $V$ of $Q$. Thus $J(u)=0$, and injectivity of $J$ gives $u=0$. This proves that $\Phi_{C,D}$ is an isomorphism. Its restriction identifies $\nplus(A)$ with $\mathfrak l(C,D)$ since it sends $e_i$ to $\theta_i$.
\end{proof}

\section{Nonrealization and Serre quotients}\label{app:serre-quotients}
Put $A=A(C,D)$ and let
\begin{equation}\label{eq:app-exceptional-edges}
    \mathscr B=
    \bigl\{
        \{i,j\}:c_{ij}<0,\ d_i,d_j>1,\
        \gcd(d_i,d_j)=1
    \bigr\}.
\end{equation}
Thus condition~\eqref{eq:edge-condition} holds if and only if $\mathscr B=\varnothing$. To prove Theorem~\ref{thm:nonrealization}, we first show that
$(\ad\theta_i)^{-a_{ij}}\theta_j\ne0$ for all $i\ne j$. We then construct a nonzero element of $\ker\Phi_{C,D}$ for each edge in $\mathscr B$.

\subsection{Proof of Theorem~\ref{thm:nonrealization}}
\begin{proposition}\label{prop:sharp-Serre}
For all distinct vertices $i,j\in I$, we have
\[
    (\ad\theta_i)^{-a_{ij}}(\theta_j)\ne0.
\]
\end{proposition}
\begin{proof}
    If $c_{ij}=0$, the assertion is immediate. If the restriction of $(C,D)$ to $\{i,j\}$ satisfies condition~\eqref{eq:edge-condition}, the assertion follows from Theorem~\ref{thm:main}. Therefore, we assume that 
    \[
    c_{ij}<0,\qquad p=d_i>1,\qquad q=d_j>1,\qquad \gcd(p,q)=1.
    \]
Put $g=g_{ij}$ and $N=|c_{ij}|=gq$. Following the construction in \cite[Lemma~6.1]{GLS18},
we construct a locally free $\Pi$-module $X$ on which $(\ad\theta_i)^{-a_{ij}}\theta_j$ is nonzero.

Take
\[
    X_i=P\oplus Q,\qquad P\cong Q\cong H_i^N,
    \qquad X_j=H_j,
\]
with $H_i$-bases $u_{\ell,t}$ of $P$ and $v_{\ell,t}$ of $Q$, where $1\le\ell\le g$ and $0\le t<q$.

Write
\[
    u_{\ell,t}^{(s)}=\varepsilon_i^s u_{\ell,t},
    \qquad
    v_{\ell,t}^{(s)}=\varepsilon_i^s v_{\ell,t}
    \qquad(0\le s<p),
\]
for the corresponding $\C$-bases of $P$ and $Q$, and let 
\[
    w_t=\varepsilon_j^t \qquad (0\le t<q)
\]
be the $\C$-basis of $X_j=H_j$.

The loops act by
\[
    \varepsilon_i u_{\ell,t}^{(s)}
    =
    \begin{cases}
        u_{\ell,t}^{(s+1)}, & 0\le s<p-1,\\
        0, & s=p-1,
    \end{cases}
    \qquad
    \varepsilon_i v_{\ell,t}^{(s)}
    =
    \begin{cases}
        v_{\ell,t}^{(s+1)}, & 0\le s<p-1,\\
        0, & s=p-1,
    \end{cases}
\]
and
\[
    \varepsilon_j w_t
    =
    \begin{cases}
        w_{t+1}, & 0\le t<q-1,\\
        0, & t=q-1.
    \end{cases}
\]
Define the arrow maps by
\[
    a_\ell(w_t)=u_{\ell,t}^{(p-1)},
    \qquad
    b_\ell(v_{\ell,t}^{(0)})=w_t
    \qquad(1\le\ell\le g,\ 0\le t<q),
\]
and all other basis vectors are mapped to zero. It follows that
\[
    \varepsilon_i^p=\varepsilon_j^q=0,\qquad
    b_\ell a_\ell=0,\qquad
    \varepsilon_i a_\ell=b_\ell\varepsilon_i=0.
\]
Thus these maps define a locally free $\Pi$-module $X$.

Since $\Ext^1_\Pi(E_i,E_i)=0$ and $\chi_c\bigl(\Fl^E_{(i,\ldots,i)}(E_i^{\oplus r})\bigr)=r!$, we have $\theta_i^{*r}=r!\,\mathbf1_{E_i^{\oplus r}}$. We now compute
\begin{align*}
    \bigl((\ad\theta_i)^{2N}\theta_j\bigr)(X)
    &=
    \sum_{s=0}^{2N}(-1)^{2N-s}\binom{2N}{s}
    \bigl(\theta_i^{*s}*\theta_j*
          \theta_i^{*(2N-s)}\bigr)(X)\\
    &=
    (2N)!\sum_{s=0}^{2N}(-1)^s
    \bigl(
        \mathbf1_{E_i^{\oplus s}}
        *\mathbf1_{E_j}
        *\mathbf1_{E_i^{\oplus(2N-s)}}
    \bigr)(X).
\end{align*}
For each $s$, we consider partial flags
\[
    0\subseteq(U,0)\subseteq(U,H_j)\subseteq X,
\]
where $U$ is a free $H_i$-direct summand of $X_i$ of rank $s$. By the definitions of $a_\ell$ and $b_\ell$, $(U,0)$ and $(U,H_j)$ are submodules of $X$ if and only if
\[
    \varepsilon_i^{p-1}P
        \subseteq U\subseteq P\oplus\varepsilon_iQ.
\]
Since $U$ is a direct summand of $X_i$, the first inclusion gives
\[
    N=\dim_\C\varepsilon_i^{p-1}P
      \le\dim_\C(U \cap \varepsilon_i^{p-1}X_i)
      =\dim_\C\varepsilon_i^{p-1}U=s.
\]
The second inclusion gives
\[
    s=\dim_\C(U/\varepsilon_iU)\le \dim_\C(P/\varepsilon_iP)=N.
\]
Hence $s=N$. 

Write $\pi_P$ and $\pi_Q$ for the projections from $P\oplus Q$ onto $P$ and $Q$, respectively. The projection $\pi_P|_U:U\to P$ induces an isomorphism
\[
    U/\varepsilon_iU \xrightarrow{\sim} P/\varepsilon_iP.
\]
By Nakayama's lemma, $\pi_P|_U$ is surjective, so it is an isomorphism.

Let $\mathcal U$ be the variety of these submodules. For $\varphi\in\Hom_{H_i}(P,\varepsilon_iQ)$, set
\[
    U_\varphi=\{x+\varphi(x)\mid x\in P\}.
\]
The map $x\mapsto x+\varphi(x)$ is an $H_i$-linear isomorphism, so $U_\varphi$ is free of rank $N$. It follows that
\[
    \varepsilon_i^{p-1}P \subseteq U_\varphi \subseteq P\oplus\varepsilon_iQ.
\]
Therefore the map
\[
    F:\Hom_{H_i}(P,\varepsilon_iQ)
    \longrightarrow\mathcal U,
    \qquad \varphi\longmapsto U_\varphi,
\]
is well-defined.

Conversely, define
\[
    G: \mathcal U
    \longrightarrow \Hom_{H_i}(P,\varepsilon_iQ),
    \qquad U\longmapsto\pi_Q\circ(\pi_P|_U)^{-1}.
\]
The maps $F$ and $G$ are mutually inverse and regular in the Grassmannian chart of complements to $Q$. Thus
\[
    \mathcal U\cong\Hom_{H_i}(P,\varepsilon_iQ)
        \cong\mathbb A^{(p-1)N^2}.
\]
Since $\chi_c(\mathbb A^{(p-1)N^2})=1$, we obtain
\[
    \bigl((\ad\theta_i)^{-a_{ij}}\theta_j\bigr)(X)
      =(-1)^N(2N)!\ne0.
\]
\end{proof}

We now show that each edge in $\mathscr B$ gives a nonzero element of $\ker\Phi_{C,D}$.

\begin{proposition}\label{prop:app-extra-relations}
Let $\{i,j\}\in\mathscr B$, and set
$a=-a_{ij}$ and $b=-a_{ji}$. Then
\[
    (\ad\theta_i)^{ab-2}(\ad\theta_j)^b\theta_i=0,
\]
whereas
\[
    (\ad e_i)^{ab-2}(\ad e_j)^b e_i\ne0
    \qquad\text{in }\mathfrak n_+(A).
\]
In particular, if condition~\eqref{eq:edge-condition}
fails, then $\Phi_{C,D}$ is not injective.
\end{proposition}

\begin{proof}
    Since $\gcd(d_i,d_j)=1$, we have $a=2g_{ij}d_j$ and $b=2g_{ij}d_i$. Put $N=b|c_{ij}|=ab/2$, and set
    \[
    u=(\ad\theta_j)^b \theta_i, \qquad z=(\ad \theta_i)^{ab-2}u.
    \]
    Since $z$ is primitive and homogeneous of rank $(ab-1)\mathbf e_i+b\mathbf e_j$, it suffices to show that $z(X)=0$ for every indecomposable locally free module $X$ of this rank. Then we have
    \[
    z=\sum_{k=0}^{ab-2}(-1)^k\binom{ab-2}{k}\theta_i^{*(ab-2-k)}*u*\theta_i^{*k}.
    \]
    For each $k$, consider a flag
    \[
    0\subseteq X'\subseteq X''\subseteq X
    \]
such that $X'\cong E_i^{\oplus(ab-2-k)}$, $X/X''\cong E_i^{\oplus k}$ and $Y=X''/X'$ is locally free of rank $\mathbf e_i+b\mathbf e_j$. We show that $E_j$ is a direct summand of every such $Y$. This will give $u(Y)=0$ by Lemma~\ref{lem:primitive-support} and hence $z(X)=0$.

Write $\alpha=X_{i,\mathrm{in}}$ and $\beta=X_{i,\mathrm{out}}$. Then
     \[
     \widetilde X_i\cong H_i^N,\qquad X_i\cong H_i^{2N-1},
    \qquad
    \alpha:H_i^N\longrightarrow H_i^{2N-1},\quad
    \beta:H_i^{2N-1}\longrightarrow H_i^N,\quad
    \alpha\beta=0.
     \]

By Smith normal form, each of $\Im \alpha$ and $\Im \beta$ can be generated by at most $N$ elements. If both could be generated by at most $N-1$ elements, then Lemma~\ref{lem:general splitting} would give $E_i$ as a direct summand of $X$, since $2N-1>2(N-1)$. This contradicts the indecomposability of $X$.

Thus at least one of $\Im \alpha$ and $\Im \beta$ is minimally generated by $N$ elements. Suppose first that $\Im \alpha$ is minimally generated by $N$ elements as an $H_i$-module.

Let $u_1,\ldots,u_N$ be an $H_i$-basis of $\widetilde X_i$. Then $\Im \alpha=\sum_{r=1}^N H_i \alpha(u_r)$. We claim that $\Ker \alpha \subseteq \varepsilon_i \widetilde X_i$. Indeed, for $v=\sum_{r=1}^N c_r u_r \in \Ker \alpha$, we have
\[
    \sum_{r=1}^N c_r\alpha(u_r)=0.
\]
If some $c_r$ were invertible, then $\alpha(u_r)$ would be generated by the remaining $N-1$ elements, contradicting the minimality of $N$. 
Thus every $c_r$ belongs to $\varepsilon_iH_i$, and hence $v\in\varepsilon_i\widetilde X_i$, as claimed. 
Since $\alpha\beta=0$, we obtain
\[
    \Im\beta\subseteq\Ker\alpha
    \subseteq\varepsilon_i\widetilde X_i.
\]
Consequently, 
\[
    \beta(\varepsilon_i^{d_i-1}X_i)
    =\varepsilon_i^{d_i-1}\beta(X_i)=0.
\]
Similarly, suppose that $\Im\beta$ is minimally generated by $N$ elements as an $H_i$-module. By the same proof we have $\Im\alpha\subseteq\varepsilon_i X_i$.

Thus $X$ satisfies at least one of
\[
    \text{(a)}\quad
    \Im X_{i,\mathrm{in}} \subseteq\varepsilon_iX_i,
    \qquad
    \text{(b)}\quad X_{i,\mathrm{out}} (\varepsilon_i^{d_i-1}X_i)=0.
\]
Now consider the flag as above and put $Y=X''/X'$. We have $X'_j=0$ and $X''_j=X_j$.

In case (a), since $X''$ is a submodule and $X''_j=X_j$, we have $\Im\alpha\subseteq X''_i$. Moreover, $X_i/X''_i$ is free over $H_i$, so
\[
    X''_i\cap\varepsilon_iX_i=\varepsilon_iX''_i.
\]
Thus $\Im\alpha\subseteq \varepsilon_iX''_i$. Passing to the quotient by $X'_i$ gives 
\[
    \Im Y_{i,\mathrm{in}}\subseteq\varepsilon_iY_i.
\]
In case (b), since $X'$ is a submodule and $X'_j=0$, we have $\beta(X'_i)=0$. Thus $\beta|_{X''_i}$ induces $Y_{i,\mathrm{out}}$ on $Y_i=X''_i/X'_i$. Since $\beta(\varepsilon_i^{d_i-1}X_i)=0$, we obtain
\[
    Y_{i,\mathrm{out}}(\varepsilon_i^{d_i-1}Y_i)=0.
\]
We show that either condition makes $Y$ decomposable. Since $Y$ is locally free of rank $\mathbf e_i+b\mathbf e_j$, we have
\[
    Y_i\cong H_i,\qquad Y_j\cong H_j^b,
    \qquad
    \widetilde Y_j
    ={}_jH_i\otimes_{H_i}Y_i
    \cong H_j^{|c_{ji}|}.
\]
Thus these maps have the form
\[
    \widetilde Y_j
    \xrightarrow{\,Y_{j,\mathrm{in}}\,}
    Y_j
    \xrightarrow{\,Y_{j,\mathrm{out}}\,}
    \widetilde Y_j.
\]
Since $f_{ji}=d_i$ and $f_{ij}=d_j$, all commutation relations $\varepsilon_j^{f_{ij}}\alpha_{ji}^{(\ell)}=\alpha_{ji}^{(\ell)}\varepsilon_i^{f_{ji}}$ hold with both sides equal to zero. Hence
\[
    {}_jH_i\cong
    (H_j\otimes_{\mathbb C}H_i)^{\oplus g_{ij}}
\]
as $(H_j,H_i)$-bimodules, and therefore
\[
    \widetilde Y_j
    ={}_jH_i\otimes_{H_i}Y_i
    \cong H_j\otimes_{\mathbb C}Y_i^{\oplus g_{ij}}.
\]

In case (a), each arrow map from $j$ to $i$ sends $Y_j$ into $\varepsilon_iY_i$. The definition of $Y_{j,\mathrm{out}}$ gives
\[
    \Im Y_{j,\mathrm{out}}
    \subseteq H_j\otimes_{\mathbb C}
    (\varepsilon_iY_i)^{\oplus g_{ij}}.
\]
The module on the right is free over $H_j$ of rank $g_{ij}(d_i-1)=|c_{ji}|-g_{ij}$. Hence $\Im Y_{j,\mathrm{out}}$ can be generated by at most $|c_{ji}|-g_{ij}$ elements as an $H_j$-module.

In case (b), each arrow map from $i$ to $j$ vanishes on $\varepsilon_i^{d_i-1}Y_i$. The definition of $Y_{j,\mathrm{in}}$ gives
\[
    H_j\otimes_{\mathbb C}
    (\varepsilon_i^{d_i-1}Y_i)^{\oplus g_{ij}}
    \subseteq\Ker Y_{j,\mathrm{in}}.
\]
The module on the left is free over $H_j$ of rank $g_{ij}$. Hence $\Im Y_{j,\mathrm{in}}$ can be generated by at most $|c_{ji}|-g_{ij}$ elements as an $H_j$-module.

In either case, one of the two images can be generated by at most $|c_{ji}|-g_{ij}$ elements, and the other by at most $|c_{ji}|$ elements, as $H_j$-modules. Since
\[
    \rank_{H_j}Y_j=b=2|c_{ji}|
    >|c_{ji}|+(|c_{ji}|-g_{ij}),
\]
Lemma~\ref{lem:general splitting} gives $E_j$ as a direct summand of $Y$. Since $Y_i\ne0$, the module $Y$ is decomposable. Thus $u(Y)=0$ and it follows that $z=0$.

Finally, to prove that $(\ad e_i)^{ab-2}(\ad e_j)^b e_i\ne0$, we apply
$(\ad f_j)^b(\ad f_i)^{ab-2}$ to this element. Since $b\ge2$, we have
\begin{align*}
    [f_i,(\ad e_j)^b e_i]&=-(\ad e_j)^b h_i=0,\\
    [h_i,(\ad e_j)^b e_i] &=(2-ab)(\ad e_j)^b e_i.
\end{align*}
Hence we have
\begin{align*}
    (\ad f_j)^b(\ad f_i)^{ab-2}(\ad e_i)^{ab-2}(\ad e_j)^b e_i&=((ab-2)!)^2
  (\ad f_j)^b(\ad e_j)^b e_i\\
&=((ab-2)!)^2(b!)^2e_i\ne0.
\end{align*}
\end{proof}

\begin{proof}[Proof of Theorem~\ref{thm:nonrealization}]
Assume that $(C,D)$ does not satisfy condition~\eqref{eq:edge-condition}. Put $L=\mathfrak l(C,D)$. Suppose that $L\cong\mathfrak n_+(B)$ for some symmetrizable generalized Cartan matrix $B$. 

Let $\{L^n\}_{n\ge1}$ denote the lower central series of $L$, where $L^1=L$ and $L^n=[L,L^{n-1}]$ for $n\ge2$. We use the same notation for $\mathfrak n_+(A)$ and $\mathfrak n_+(B)$. By Proposition~\ref{prop:app-extra-relations}, the canonical surjection
\[
    \pi:\mathfrak n_+(A)\longrightarrow L,
    \qquad e_i\longmapsto\theta_i,
\]
is not injective. Choose a nonzero homogeneous element $z\in\ker\pi$, and let $d_z$ be its total degree. Let $N=\max\left\{ d_z,\ \max_{i\ne j}(1-a_{ij}),\ \max_{i\ne j}(1-b_{ij})\right\}$ and set 
\[
L_N=L/L^{N+1}, \qquad A_N=\nplus(A)/\nplus(A)^{N+1} ,\qquad B_N=\nplus(B)/\nplus(B)^{N+1}.
\]
These are finite-dimensional nilpotent Lie algebras, and the assumed isomorphism induces $L_N\cong B_N$.

Write $\bar\theta_i$ for the image of $\theta_i$ in $L_N$. The rank grading of $L$ descends to $L_N$ and defines commuting semisimple derivations $D_i$, $i\in I$, satisfying
\[
    D_i(\bar\theta_j)=\delta_{ij}\bar\theta_j.
\]
Thus
\[
    T=\bigoplus_{i\in I}\C D_i
\]
is a torus of dimension $|I|$. Since $\dim_\C L_N/[L_N,L_N]=|I|$, it follows from
\cite[Lemma~2.8]{San82} and \cite[Definition~2.9]{San82} that $T$ is maximal and $L_N$ is of maximal rank. $(\bar\theta_i)_{i\in I}$ is a $T$-msg in the sense of \cite[Definition~2.5]{San82}.

By Proposition~\ref{prop:A-Serre} and Proposition~\ref{prop:sharp-Serre}, we have
\[
    (\ad\bar\theta_i)^{-a_{ij}}\bar\theta_j\ne0,
    \qquad
    (\ad\bar\theta_i)^{1-a_{ij}}\bar\theta_j=0
    \qquad(i\ne j).
\]
Thus $A$ is the  Cartan matrix associated with $L_N$ in the sense of \cite[Theorem~3.2]{San82}. The same argument applied to $\nplus(B)$ shows that $B_N$ is of maximal rank with associated Cartan matrix $B$. Since $L_N\cong B_N$, \cite[Theorem~3.2]{San82} implies that $A$ and $B$ differ by a simultaneous permutation.

The map $\pi$ induces a surjection $\bar\pi: A_N \twoheadrightarrow L_N$. Let $\bar z$ be the image of $z$ in $A_N$. We have $\bar z\ne0$, while $\bar\pi(\bar z)=0$. Thus $\bar\pi$ is not injective. Therefore
\[
    \dim_\C L_N< \dim_\C A_N =\dim_\C B_N,
\]
contrary to $L_N\cong B_N$. This proves the assertion for $L$.

Lemma~\ref{lem:enveloping-realization} gives $\mathcal M(C,D)\cong U(L)$ as Hopf algebras. Suppose that $\mathcal M(C,D)\cong U(\mathfrak n_+(B))$ as unital $\C$-algebras for some symmetrizable generalized Cartan matrix $B$. Then $U(L)\cong U(\mathfrak n_+(B))$ as unital $\C$-algebras, and \cite[Proposition~4.1]{RU07} implies that
\[
    \operatorname{gr}L
    \cong\operatorname{gr}\mathfrak n_+(B),
\]
where $\operatorname{gr}L=\bigoplus_{r\ge1}L^r/L^{r+1}$. Since $L$ and $\nplus(B)$ are positively graded and generated in degree one, their total gradings give Lie algebra isomorphisms
\[
    L\cong\operatorname{gr}L,\qquad
    \mathfrak n_+(B)\cong\operatorname{gr}\mathfrak n_+(B).
\]
Hence $L\cong\mathfrak n_+(B)$, contradicting the assertion for $L$ proved above.
\end{proof}

\subsection{Serre quotient realizations}
We now return to an arbitrary pair $(C,D)$. Put $\mathcal M=\mathcal M(C,D)$ and let $\mathscr B$ be as in~\eqref{eq:app-exceptional-edges}.
Define $B=B(C,D)=(b_{ij})_{i,j\in I}$ by
\begin{equation}\label{eq:app-matrix-B}
 b_{ii}=2,\qquad
 b_{ij}=
 \begin{cases}
 c_{ij},&\{i,j\}\in\mathscr B,\\
 a_{ij},&\{i,j\}\notin\mathscr B.
 \end{cases}
\end{equation}
The positive integral diagonal matrix
\[
 D_A=\diag(2d_i/\lambda_i)_{i\in I}
\]
symmetrizes both $A$ and $B$. In particular, $B$ is a symmetrizable generalized Cartan matrix.

Let $\mathcal J_{\mathscr B}$ be the two-sided ideal of $\mathcal M$ generated by
\begin{equation}\label{eq:app-exceptional-ideal}
 (\ad\theta_i)^{1-c_{ij}}\theta_j
\end{equation}
for all ordered pairs $(i,j)$ with $\{i,j\}\in\mathscr B$. Since these elements are primitive, $\mathcal J_{\mathscr B}$ is a Hopf ideal.

\begin{lemma}[An intermediate Serre quotient] \label{lem:app-intermediate-quotient}
For every pair $(C,D)$, there is a canonical rank-graded Hopf algebra isomorphism
\begin{equation}\label{eq:app-intermediate-isomorphism}
\begin{aligned}
 U(\nplus(B))&\xrightarrow{\sim}
\mathcal M(C,D)/\mathcal J_{\mathscr B},\\
 e_i&\longmapsto\theta_i+\mathcal J_{\mathscr B}.
\end{aligned}
\end{equation}
\end{lemma}

\begin{proof}
    The $B$-Serre relations follow from Proposition~\ref{prop:A-Serre} and the definition of $\mathcal J_{\mathscr B}$. Thus the assignment in the statement induces a rank-graded Hopf algebra homomorphism
\[
    \psi_B:U(\nplus(B))
    \longrightarrow
    \mathcal M(C,D)/\mathcal J_{\mathscr B},
    \qquad
    e_i\longmapsto\theta_i+\mathcal J_{\mathscr B}.
\]
Since the classes $\theta_i+\mathcal J_{\mathscr B}$ generate the quotient as a unital algebra, $\psi_B$ is surjective.

For injectivity, we modify the auxiliary quiver $Q$ of Section~\ref{sec:injectivity}. Keep its vertex set and all arrows on edges outside $\mathscr B$. For each $\{i,j\}\in\mathscr B$, choose a direction $j\to i$ and retain only the arrows
\[
    \alpha_{ij,w}^{(\ell)}:
    (j,[w]_{m_j})\longrightarrow(i,[w]_{m_i}),
    \qquad
    0\le w<L,\quad 1\le\ell\le g_{ij}.
\]
Omit the arrows $\beta_{ij,w}^{(\ell)}$. Denote the resulting quiver by $Q'$ and its symmetric generalized Cartan matrix by $\widetilde B$.

For $\{i,j\}\in\mathscr B$, we have $f_{ij}=d_j$ and $f_{ji}=d_i$. For each fixed vertex $(j,b)$, it follows that
\[
    \sum_{a=0}^{m_i-1}
        \widetilde b_{(i,a),(j,b)}
    =-g_{ij}\frac{L}{m_j}
    =-g_{ij}d_j
    =c_{ij}
    =b_{ij}.
\]
Interchanging $i$ and $j$ gives the corresponding identity for $b_{ji}$. For all other edges, the arrows are unchanged, so the folding identities follow from Section~\ref{sec:injectivity}.

$Q'$ satisfies the hypotheses of Lemma~\ref{lem:inj-cover}. Combining the resulting embedding with Theorem~\ref{thm:inj-hall}, we obtain an injective algebra homomorphism
\[
\begin{aligned}
    J_B:U(\nplus(B))&\hookrightarrow\mathcal M(Q'),\\
    e_i&\longmapsto
    \sigma_i:=
    \sum_{a=0}^{m_i-1}\mathbf1_{S_{(i,a)}}.
\end{aligned}
\]
Regard a representation $V$ of $Q'$ as a representation of $Q$ by assigning zero maps to the deleted arrows. Define $F(V)$ using the formulas of Section~\ref{sec:injectivity}. For $\{i,j\}\in\mathscr B$, we have $k_{ij}=1$, so $t=L$, the weight shift is zero, and $b_\ell=0$. Since $f_{ji}=d_i$ and $f_{ij}=d_j$, the commutation relations on this edge hold by loop nilpotence, and its mesh contributions vanish. For all other edges, the verification in Lemma~\ref{lem:relations} applies unchanged. Thus $F(V)$ is a locally free $\Pi(C,D)$-module, and all its operators remain homogeneous.

By the same argument as in the proof of Lemma~\ref{lem:flag-comparison}
(see also \cite[Appendix~A.2]{Fei26}), we obtain
\[
    (\theta_{i_1}*\cdots*\theta_{i_r})(F(V))
    =
    (\sigma_{i_1}*\cdots*\sigma_{i_r})(V)
\]
for every word and every representation $V$ of $Q'$. Consequently, the assignment $\theta_i\mapsto\sigma_i$ induces an algebra homomorphism
\[
    \rho:\mathcal M(C,D)\longrightarrow\mathcal M(Q').
\]
Since the $\sigma_i$ satisfy the $B$-Serre relations, $\rho$ annihilates $\mathcal J_{\mathscr B}$ and induces a homomorphism
\[
    \bar\rho:
    \mathcal M(C,D)/\mathcal J_{\mathscr B}
    \longrightarrow\mathcal M(Q').
\]
By construction, $\bar\rho\circ\psi_B=J_B$. Since $J_B$ is injective, $\psi_B$ is injective. This proves the claimed isomorphism.
\end{proof}

As a consequence, we recover the full $C$-Serre quotient realization established by Fei \cite[Theorem~10.1]{Fei26}. In particular, this gives the algebra isomorphism in \cite[Theorem~7.4(i)]{GLS18} without assuming the Support Conjecture.

\begin{corollary}[The full $C$-Serre quotient] \label{cor:app-full-Serre-quotient}
Let $\mathcal I_C$ be the two-sided ideal of $\mathcal M(C,D)$ generated by
\[
    (\ad\theta_i)^{1-c_{ij}}\theta_j,
    \qquad i\ne j.
\]
For every pair $(C,D)$, there is a canonical rank-graded Hopf algebra isomorphism
\[
\begin{aligned}
    U(\nplus(C))
    &\xrightarrow{\sim}
    \mathcal M(C,D)/\mathcal I_C,\\
    e_i&\longmapsto\theta_i+\mathcal I_C.
\end{aligned}
\]
\end{corollary}

\begin{proof}
Since $\mathcal J_{\mathscr B}\subseteq\mathcal I_C$ and $b_{ij}\le c_{ij}$ for $i\ne j$, the assertion follows from Lemma~\ref{lem:app-intermediate-quotient} and the Serre presentation of $U(\nplus(C))$.
\end{proof}

\section*{Declaration of AI-assisted research and manuscript preparation}
The author used ChatGPT for literature searches, language editing, computations in examples, and discussions of proof strategies. This assistance helped identify a counterexample to an earlier formulation of Theorem~\ref{thm:main}. Further computations led to Proposition~\ref{prop:app-extra-relations}, which the author used to prove Theorem~\ref{thm:nonrealization}. The author has carefully checked the paper and takes full responsibility for its correctness. 

\end{document}